# Primitive representations of finite semigroups 1

Steve Seif  Johnny Ray Sena

August 11, 1997

## Abstract

Primitive representations of finite groups as well as primitive finite groups were classified in the O'nan-Scott Theorem. In this paper we classify faithful finite primitive semigroup representations. To each finite primitive representation, we associate an invariant, a finite dimensional matrix with entries in a primitive finite group representation. To a large extent, this matrix determines the representation. In later papers in this series the invariant is further explored.

Our primitivity results rely on two main ideas, one of which is a small but important part of the theory of tame algebras developed by R. McKenzie, while the other is our adaptation of Rees matrix theory to the representation theory of finite semigroups. Using the latter idea we are able to provide a description of all representations of a regular $\mathcal{J}$ class of a finite semigroup.

## 1  Introduction

We present a surprisingly cohesive picture of the finite primitive semigroup representations. This paper is the first in a series of papers on the subject. As a byproduct of the work on primitivity, we are also able to make a contribution to the general problem of describing representations of finite semigroups by semigroups of transformations.

We begin by reviewing the main notions associated with representing a (subset of a) semigroup by a (subset of a) semigroup of transformations on a set. Let $Y$ be a set. Recall that $T_Y$ denotes the semigroup of transformations of $Y$ and that $S_Y$ denotes the subgroup of $T_Y$ consisting of the permutations of $Y$. We let $Y$ do double duty: it denotes both the set $Y$ and the set of constant functions on the set $Y$. For a semigroup $S$ and a subset $U$ of $S$, by a representation of $U$ we mean a map $\gamma : U \to T_Y$ which satisfies the following. For $s, t \in U$, if $st \in U$ then $\gamma(st) = \gamma(s)\,\gamma(t)$. So $\gamma$ is a partial homomorphism from the partial semigroup $U$ into the semigroup $T_Y$. Given $U$ a subset of a semigroup $S$ with a representation $\gamma : U \to S_Y$, our notation is $(Y; U; \gamma)$; if there is only one representation of $U$ under discussion, then we often abbreviate $(Y; U; \gamma)$ to $(Y; U)$. If $U$ is itself a semigroup, then we say that the representation is called a *semigroup representation*.

Let $Y$ be a finite set. Given a subset $W$ of $T_Y$, we say that $(Y;W)$ is an **action**. Of course any action can be turned into a representation by considering $W$ as a subset of $\langle W \rangle$, the semigroup of transformations of $Y$ generated by $W$. But a great many of the results of this paper and its sequels begin with an abstract semigroup and and present a description of certain representations of its subsets.

A representation $(Y;U;\gamma)$ is *faithful* if $\gamma$ is one-to-one. For a representation $(Y;U)$ an equivalence relation $\theta$ of $Y$ is said to be *U-compatible* if for all $(a,b) \in \theta$ and all $u \in U$, $(u(a), u(b)) \in \theta$. A representation $(Y;U)$ is said to be *primitive* if the only U-compatible equivalence relations are the diagonal and the universal relations. Observe that if $(Y;S)$ is a representation and $|Y| = 1$ or $2$, then $(Y;S)$ is primitive. Thus our main theorem, Theorem 1.8, concerns representations on sets with more than 2 elements. Theorem 1.8 involves a finite primitive group representation. These were classified in the O'nan-Scott Theorem. See [5] for a proof of the O'nan-Scott Theorem.

Let $(Y;U)$ be a representation. The *expansion of U by constants* (denoted by $U^Y$) is the representation consisting of the transformations in $U$ along with the full set of constant transformations of $Y$. Note that $(Y;U)$ is primitive if and only if $(Y;U^Y)$ is primitive.

The introduction continues with the definitions and notation involved in the main theorem of this paper, Theorem 1.8, whose statement essentially ends the introduction. In Section 2, we provide some background on semigroup theory. (see [1] for additional background). In Section 3 we introduce the R-representations, the representations of regular $\mathcal{J}$ classes of finite semigroups and prove a Fundamental Lemma for R-representations. In Section 4 we prove our results on primitive semigroup representations. Ideas from the theory of tame algebras (see [3] and [6]) are the springboard for many of the results in the section. In Section 5, we explore what we call the **c-ramified matrices** associated with a finite primitive representation. In Section 6, we return to R-representations and classify these completely. Some general remarks and problems are offered in Section 7.

## 1.1 Primitive semigroup representations

Our description of primitive representations relies upon the *c-ramified Rees representations* which we describe in Definitions 1.1 and 1.2.

**Definition 1.1.** *Let $(X;G)$ be a group representation and let $m,n$ be positive integers. An $m \times n$ matrix $P$ with entries in $X \cup G$ is said to be a matrix which is **c-ramified over** $(X;G)$. We say that $P$ is **regular** if every row and every column of $P$ has at least one entry from $G$.*

*Entries of $P$ are thought of as functions on the set $X$; entries from the set $X$ are constant functions, while entries from $G$ are permutations.*

We depart from c-ramified matrices for the moment. Whenever we refer to an $m \times n$ matrix $A$, we let $\Lambda = \{1, \ldots, \lambda, \ldots, m\}$ index the rows of $A$; we let $I = \{1, \ldots, i, \ldots, n\}$ index its columns. The $(\lambda, i)$ entry of $A$ will be denoted by $A(\lambda, i)$. The $\lambda$ row of $A$ will be denoted by $A_\lambda$; the i column of $A$ will be denoted by $A^{(i)}$. Recall that a matrix with exactly one non-0 entry (in each row) (in each column) is called a *monomial (row monomial) (column monomial)* matrix. We refer to a $n \times 1$ matrix as an *n*-vector. Letting $(X; G)$ once again denote a group representation, observe that if $M$ is a $n \times m$ monomial matrix with entries in $X \cup G$ and $P$ is an $m \times n$ c-ramified matrix over $(X; G)$, then the matrix multiplication $MP$ is well-defined, under the following conventions. For any $x \in X$, $x + 0 = x = 0 + x$, and $x0 = 0 = 0x$; for $g \in G$, $g0 = 0 = 0g$. If $u \in G \cup X$ and $v \in G \cup X$, then $uv$ is defined to be the composition of $u$ with $v$, a function of $X$. Moreover, with $M$ and $P$ as above, if $v$ is a monomial $n$ vector over $X$, then under the same matrix multiplication conventions, $MPv$ is a well-defined $n$ vector over $X$. For the remainder of the paper, matrix products will be defined in terms of the conventions above. Note also that if $M, N \in J((X; G); P)$, then $MPN$ is a well defined $m \times n$ monomial matrix with a unique non-0 entry in $X \cup G$.

**Definition 1.2.** Let $(X; G)$ be a group representation and let $P$ be a regular c-ramified $m \times n$ matrix over $(X; G)$.

1. $J((X; G); P)$ will refer to the set of $n \times m$ monomial matrices over $G$. That is, each matrix $M \in J((X; G); P)$ contains a permutation from $G$ in one entry and all the other entries of $M$ are 0.

2. Let $V$ be the set of monomial *n*-vectors over $X$.

3. Let $J((X; G); P)$ act on $V$ as follows. For $M \in J((X; G); P)$ and $v \in V$, let $M(v) = MPv$.

**Notation 1.3.**  1. For $M \in J((X; G); P)$, it is often convenient to use a 3-tuple form to describe $M$. If the non-0 entry of $M$ is $g \in G$ and $g$ occurs in the $(\lambda, r)$ entry, then the 3-tuple form of $M$ will be $[r, g, \lambda]$.

2. We use $r_x$ ($1 \leq r \leq n$, $x \in X$) to denote the element of $V$ with a non-0 entry $x$ in the $r$ row.

3. The action of $J((X; G); P)$ on $V$ can be described in terms of the notation of the two previous items. For $[r, g, \lambda] \in J((X; G); P)$ and $i_x \in V$, $[r, g, \lambda](i_x) = r_{gP(\lambda, i)(x)}$.

We often use $i, r, s$ as elements of $I$, $\alpha, \beta, \lambda$ as elements of $\Lambda$, and $a, b, x, y$ as elements of $X$. We use the matrix and 3-tuple form for elements of $J((X; G); P)$ interchangably.

**Definition 1.4.** Let $P$ be a $n \times m$ ramified Rees matrix over $(X; G)$. Let $\theta_P = \{(r_x, s_y) \in V^2 : \forall \lambda \in \Lambda \ P(\lambda, r)(x) = P(\lambda, s)(y)\}$.

Observe that $\theta_P$ is $J((X; G); P)$-compatible. In fact if $(r_x, s_y) \in \theta_P$, then for all $M \in J((X; G); P)$, $M$ acts on $r_x$ and $s_y$ in the same way: $MP(r_x) = MP(s_y)$).

We abbreviate $V/\theta_P$ to $V_{\theta_P}$. By the previous paragraph, the action of $J((X;G);P)$ on $V$ determines an action of $J((X;G);P)$ on $V/\theta_P$. We denote this action by $(V_{\theta_P}; J((X;G);P))$ and regard it as an action by matrix multiplication. For $M \in J((X;G);P)$, $r_x/\theta_P \in V_{\theta_P}$, the well-defined action is $r_x/\theta_P \to (MPr_x)/\theta_P$.

The actions of $J((X;G);P)$ on $V$ and on $V_{\theta_P}$ are not generally closed under composition. If $M, N \in J((X;G);P)$, then as we've observed in the paragraph preceding Definition 1.2, $MPN$ is a monomial $n \times n$ matrix with non-0 entry in $X \cup G$. Thus the composition of the actions of $M$ and $N$ on $V$, or on $V_{\theta_P}$, can be represented by a matrix in $J((X;G);P)$ if and only if $MPN \in J((X;G);P)$ (i.e. $MPN$ has its unique non-0 element in $G$); the composition of the actions associated with $M$ and $N$ results in a constant function on $V$, or on $V_{\theta_P}$, if $MPN$ has an element of $X$ as its unique non-0 entry. If the non-0 entry of $MPN$ is $x \in X$ and $x$ occurs in the i row of $MPN$, then the image of the composition is $\{i_x\}$ in the case of $V$ and $\{i_x/\theta_P\}$ in the case of $V_{\theta_P}$.

Consider the actions $(V; J((X;G);P)^V)$ and $(V_{\theta_P}; J((X;G);P)^{V_{\theta_P}})$, which result by adding the constant transformations of $V$, $V_{\theta_P}$, respectively to $J((X;G);P)$. From the remarks of the preceding paragraph, we have that $J((X;G);P)^V$ and $J((X;G);P)^{V_{\theta_P}}$ are semigroups of transformations.

**Notation 1.5.**   1. $L((X;G);P) = J((X;G);P)^{V_{\theta_P}}$

   2. $U((X;G);P) = J((X;G);P)^V$

To a c-ramified matrix $P$, we associate a graph $Gr(P)$.

**Definition 1.6.** Let $P$ be an $m \times n$ c-ramified matrix over a group representation $(X;G)$. We define $Gr(P)$, a graph with vertex set consisting of the columns of $P$, $\{P^{(i)} : i \in I\}$. For $1 \leq u \leq n$, $r \neq s \in I$, we let $P^{(r)}, P^{(s)}$ be an edge of $Gr(P)$ if there exits $x, y \in X$ such that for all $\lambda \in \Lambda$, $P(\lambda, r)(x) = P(\lambda, s)(y)$.

In particular, if $r \neq s \in I$, $P^{(r)}, P^{(s)}$ is an edge of $Gr(P)$ if and only if there exists $x, y \in X$ such that $(r_x, s_y) \in \theta_P$.

**Definition 1.7.** Let $P$ be an $m \times n$ c-ramified matrix over the group representation $(X;G)$.

- For $\alpha \neq \beta \in \Lambda$, we say that $M_\alpha$ and $M_\beta$ are **proportional** if there exists $g \in G$ such that for for all $i \in I$, $gP(\alpha, i) = P(\beta, i)$.

- For $r \neq s \in I$, we say that $M^{(r)}$ and $M^{(s)}$ are proportional if there exists $h \in G$ such that for all $\lambda \in \Lambda$, $P(\lambda, r)h = P(\lambda, s)$.

- If $P$ has no proportional rows, then we say that $P$ is **right reductive**. If $P$ has no proportional columns, then we say that $P$ is **left reductive**. If $P$ is left and right reductive, the we say that $P$ is **reductive**.

Given a representation $(Y; S; \gamma)$ as we mentioned, we often drop $\gamma$ and write $(Y; S)$. For $s \in S$ and $y \in Y$, we replace $\gamma(s)(y)$ by $s(y)$.

We provide a (standard) definition of *equivalence* of representation. Two representations $(Y; U; \gamma) = (Y; U)$ and $(Z; V; \mu) = (Z; V)$ are said to be *equivalent* if there exists a bijection $\beta : (U \cup Y) \mapsto (V \cup Z)$ satisfying the following. For all $u \in U$ and all $y \in Y$, $\beta(u(y)) = \beta(u)\,(\beta(y))$ and for all $u, v \in U$ such that $uv \in S$, $\beta(uv) = \beta(u)\,\beta(v)$. Observe that if $(Y; U)$ and $(Z; V)$ are equivalent by a map $\beta$, then $\beta(Y) = Z$ and $\beta(U) = V$. Note that if we had used the more proper notation, the definition of equivalence of representations becomes cluttered with symbols.

Suppose that $(Y; S; \gamma) = (Y; S)$ is a faithful semigroup action. Let $\Omega(S) = \{f \in Y^Y : fS \cup Sf \subset S\}$ Usually $\Omega(S)$ denotes the so-called *translational hull* of $S$. It turns out that the two usages of $\Omega(S)$ coincide. In a companion paper [7], translational hulls of $L((X; G); P)$ and $U((X; G); P)$ semigroups are given detailed description.

Let $(Y; S)$ be a representation. A subset $W \subseteq Y$ will be said to be a *minimal set* of $Y$ if $|W| \geq 2$, there exists $s \in S$ such that $s(Y) = W$, and $K \subseteq W$, $|K| \geq 2$, and there exists $t \in S$ such that $t(Y) = K$, imply that $K = W$. So $W$ is minimal if it is minimal in the set of non-singleton image sets of $(Y; S)$. Observe that if $|Y| > 2$, then $(Y; S)$ has minimal sets.

**Theorem 1.8.** *Let $(Y; S)$ be a finite semigroup representation with $|Y| > 2$. Then $(Y; S)$ is a faithful primitive representation if and only if $(Y; S^Y)$ is equivalent to $(V_{\theta_P}; T)$, where*

1. $T$ is a subsemigroup of $\Omega(L((X; G); P)$

2. $J((X; G); P) \subseteq T$

3. $(X; G)$ is a primitive faithful group representation.

4. $P$ is a regular, reductive c-ramified matrix over $(X; G)$.

5. The graph $Gr(P)$ is connected.

*Moreover, $J((X; G); P)$ is the set of minimal functions of $(V_{\theta_P}; T)$ and for all $s \neq t \in T$, there exists $j, k \in J((X; G); P)$ such that $sj \neq tj$ and $ks \neq kt$.*

Theorem 1.8 assigns to each faithful primitive semigroup representation $(Y; S)$ a reductive regular c-ramified matrix $P$ over a faithful primitive group representation. As would be expected, $P$ is not unique, but the next proposition clarifies the situation.

Let $P$ be an $m \times n$ c-ramified matrix over a group representation $(X; G)$. Suppose that $(X'; G')$ is a group representation and $\beta : (X; G) \to (X'; G')$ is an equivalence. We let $\beta(P)$ denote the $m \times n$ matrix such that $\beta(P)(\lambda, i) = \beta(P(\lambda, i))$ ($\lambda \in \Lambda$, $i \in I$).

Let $(X; G)$ be a group representation. A $k \times k$ matrix with entries in $G^0$ which is both row and column monomial is said to be a *permutational* matrix over $G$.

**Proposition 1.9.** *Suppose that $P$ is a regular, reductive c-ramified matrix over the faithful group representation $(X; G)$ and $P'$ is a c-ramified matrix over the faithful group representation $(X'; G')$. Then the following are equivalent.*

1. *$L((X; G); P)$ is isomorphic to $L((X'; G'); P')$*

2. *$(V_{\theta_P}; L((X; G); P))$ is equivalent to $(V_{\theta_{P'}}; L((X'; G'); P'))$*

3. *There exist an $m \times m$ permutational matrix $A$ over $G$, an $n \times n$ permutational matrix $B$ over $G$, and an equivalence $\beta: (X; G) \to (X'; G')$ such that $P' = \beta(APB)$.*

Proposition 1.9 indicates that to every primitive representation is associated an invariant, an easily understood set of c-ramified Rees matrices. The c-ramified matrix invariant will be treated in further detail in [7] and [8].

This completes our description of the primitive semigroup representatives in this paper. In [7], $\Omega(L((X; G); P))$ is represented as a set of matrices acting on $V_{\theta_P}$ by matrix multiplication.

## 2 Semigroups

We review some definitions from semigroup theory. All semigroups are assumed to be finite. A semigroup $S$ has a *zero* if there exists an element $0 \in S$ such that $0s = s0 = 0$ for all $s \in S$. A semigroup $S$ has a *identity* if there exist an element $1 \in S$ such that $1s = s1 = s$ for all $s \in S$. If a semigroup has a either a zero or an identity, then in either case that element is unique. By $S^0$ we mean the minimal extension of $S$ to a semigroup with a zero; $S^1$ is the minimal extension of $S$ to a semigroup with an identity. If $T$ is a subset of $S$ and $T$ is closed under the semigroup operations, then $T$ is said to be a *subsemigroup* of $S$ (denoted by $T \leq S$). If $T$ is a subsemigroup of $S$ and $T$ is a group, then $T$ is said to be a *subgroup* of $S$. If $U$ is a subset of $S$, then the subsemigroup of $S$ generated by $U$ is denoted by $<U>$.

A subset $I \subset S$ is said to be a *ideal* of $S$ if for all $a \in S$, $aI \cup Ia \subset I$. We say that $S$ is *simple* if $S$ has exactly one non-empty ideal, namely $S$ itself; $S$ is 0-simple if and only if the only proper non-empty ideal of $S$ is $\{0\}$. A semigroup that is either simple or 0-simple will be referred to as a *(0)-simple* semigroup. A subset $K$ of $S$ is said to be a (left) (right) ideal if for all $s \in S$ $(Ks \subseteq K)$ $(sK \subseteq K)$ $Ks \cup sK \subseteq K$. Observe that if $Y$ is a set and $S \leq T_Y$, then the set of constant maps contained in $S$ is an ideal of $S$.

If $I$ is an ideal of $S$, then by $S/I$ we mean the semigroup with 0 and underlying set $(S \setminus I) \cup \{0\}$ with the following multiplication: if $s, t \in S \setminus I$ and $st \notin I$, then the product of $s$ and $t$ in $S/I$ is $st$; all other products are 0. We define *quotient semigroups* in a more general setting. An equivalence relation $\nu$ on $S$ is said to be a *congruence* of $S$ if for all $(s, t) \in \nu$ and all $u \in S$, $(us, ut) \in \nu$ and $(su, tu) \in \nu$. If $\nu$ is a congruence of $S$, then the multiplication operation on $S$ induces a well-defined associative multiplication on $U/\nu$.

We define a notion that bridges abstract semigroups and semigroup representations. Let $S$ be a semigroup. For a fixed $s \in S$ and all $t \in S$, $\rho(s)(t) = st$ and $\lambda(s)(t) = ts$. So $\rho(s)$ and $\lambda(s)$ define functions from $S$ to itself and it is easy to see that $(S; S; \rho)$ and $(S; S; \lambda)$ are both representations of $S$. We denote $\rho(s)$ by $\rho_s$ and $\lambda(s)$ by $\lambda_s$.

## 2.1 Green's Relations, completely (0)-simple semigroups

We review some basic notions concerning *Green's Relations* from the theory of finite semigroups. The sets $aS^1$, $S^1a$ and $S^1aS^1$ are called *the principal right ideal generated by a*, *the principal left ideal generated by a* and the *the principal ideal generated by a*, respectively. The following is a list of Green's Relations for finite semigroups.

1. $a\mathcal{J}b$, if $a$ and $b$ generate the same principal ideal.

2. $a\mathcal{R}b$, if $a$ and $b$ generate the same principal right ideal

3. $a\mathcal{L}b$, if $a$ and $b$ generate the same principal left ideal.

4. $a\mathcal{H}b$, if $a\mathcal{L}b$ and $a\mathcal{R}b$.

Let $a \in S$. Then the $(\mathcal{R})$ $(\mathcal{L})$ $(\mathcal{H})$ $\mathcal{J}$ class of $a$ is denoted by $(R_a)$ $(L_a)$ $(H_a)$ $J_a$.

For $a, b \in S$, we write $a \leq_J b$ if $S^1aS^1 \subset S^1bS^1$. Observe that $\leq_J$ is a pre-ordering of $S$ and that the equivalence relation associated with $\leq_J$ is of course $\mathcal{J}$. In particular, the $\mathcal{J}$ classes of $S$ admit a partial ordering.

Let $a \in S$. We let $\{b : b \leq a\}$ be denoted by $S_{a]}$ and let $\{c : c < a\}$ be denoted by $S_{a)}$. Note that $S_{a]}$ and $S_{a)}$ are both ideals of $S$ and that $S_{a]} \setminus S_{a)} = J_a$.

**Observation 2.1.** *Let $S$ be a semigroup with $s, t \in S$ and let $(Y; S; \gamma)$ be a representation of $S$ with $s, t \in S$.*

1. *If $s \mathcal{R} t$, then $im(\gamma(s)) = im(\gamma(t))$.*

2. *If $s \mathcal{L} t$, then $ker(\gamma(s)) = ker(\gamma(t))$.*

As we will see if $(Y; J; \gamma)$ is what we call a faithful R-representation, then the converses to the items of Observation 2.1 hold. The next lemma is a starting point for much of the algebraic theory of semigroups. We provide a proof.

**Lemma 2.2.** *(Green's Lemma For Finite Semigroups) Let $S$ be a finite semigroup with $a, u \in S$ and suppose that $ua \in J_a$.*

1. *Then $ua \in L_a$ and $\lambda_u$ maps $R_a$ bijectively onto $R_{ua}$. Moreover if $u \in J_a$, then $R_{ua} = R_u$.*

2. *Each $\mathcal{H}$ class contained in $R_a$ is mapped by $\lambda_u$ bijectively onto an $\mathcal{H}$ class contained in $R_b$.*

*The dual ($au \in J_a$ implies that $au \in R_a$, etc.) also holds.*

*Proof.* Let $a, u \in S$ and assume that $ua \in J_a$. By definition of the $\mathcal{J}$ relation, there exist $s, t \in S^1$ such that $suat = a$. Thus for all Positive integers $n$, we have $(su)^n at^n = a$. By finiteness we can choose $n$ so that $t^n$ is idempotent. It follows that $at^n = a$; hence, $(su)^n a = a$ which we rewrite as $(su)^{n-1} s(ua) = a$. It follows that $L_{au} = L_a$ thereby that $ua \in L_a$. With $a, u, ua$ as above, we have shown that there exists $v \in S$ such that $vua = a$. Let $b \in R_a$. There exists $w \in S^1$ such that $b = aw$. Hence, $vub = b$. The inverse of $\lambda_u|_{R_a}$ is $\lambda_v|_{R_{ua}}$. This completes the proof of the first sentence of item 1.

Assume now that $u \in J_a$. Then there exists $c, d \in S^1$ such that $cad = u$. Thus $suat = a$ implies that $scadat = a$. As in the proof of the first line of item 1, there exists a positive integer $n$ such that $a(dat)^n = a$, from which we have $c(a(dat)^n)d = u$. We have $(cad)aw = u$, where $w \in S$. Thus $ua(w) = u$ and we have shown that $R_{ua} = R_u$.

We prove item 2. If $c \in R_a$, then by item 1 we have $cy \in L_c$. Thus, $H_c = L_c \cap R_c = L_c \cap R_a$ is mapped by $\lambda_u$ into $L_c \cap R_{ua} = H_{ua}$. The fact that the various translations act bijectively implies that $H_c$ is mapped onto $H_{ua}$.

The dual formulation of the lemma is proved by dualizing the proofs above. □

Notice that one consequence of the lemma above is that a $\mathcal{J}$ class is equipartitioned by both $\mathcal{L}$ and $\mathcal{R}$ (hence also by $\mathcal{H}$).

Here are some direct consequences of Green's Lemma.

**Lemma 2.3.** *Let $J$ be a $\mathcal{J}$ class of $S$ and let $s, t \in J$ and let $c, d \in J^1$.*

1. *If $sJt \cap J \neq \emptyset$, then $sJt \cap J$ is an $\mathcal{H}$-class.*

2. *Suppose $H$ is an $\mathcal{H}$ class and $H \subset J$. If $cHd \cap J \neq \emptyset$, then $cHd \subset J$ and $cHd$ is an $\mathcal{H}$-class. In particular the map $x \mapsto cxd$ is a bijection between the $\mathcal{H}$ classes $H$ and $cHd$.*

The following definitions involve the notion of *regularity* in semigroups. Suppose that $a \in S$ and that $aba = a$. Then $b$ is said to be an *inverse* of $a$ and $a$ is said to be a *regular* element of $S$. Let $J$ ($H$) be a $\mathcal{J}$ ($\mathcal{H}$) class. If every element of $J$ ($H$) is regular, then $J$ ($H$) is said to be a *regular $\mathcal{J}$ ($\mathcal{H}$) class*. If every $\mathcal{J}$ class of $S$ is regular, then $S$ itself is said to be regular. A semigroup $S$ is said to be *completely (0)-simple* if $S$ is (0)-simple and $S$ is regular. A completely (0)-simple semigroup $S$ contains at most two $\mathcal{J}$ classes, the $\mathcal{J}$ class consisting of $\{0\}$ (if $S$ has a 0) and the regular $\mathcal{J}$ class consisting of $S \setminus \{0\}$. We refer to non-0 $\mathcal{J}$ class of a completely (0)-simple semigroup as the maximal $\mathcal{J}$ class of $S$ and usually denote it by $J$.

The following well-known results are left as a not completely trivial exercise.

**Lemma 2.4.** *Let $a \in S$.*

1. *If $a$ is regular, then all elements in $J_a$ are also regular.*

2. The following are equivalent
    (a) $J_a$ is regular
    (b) $J_a$ contains an idempotent.
    (c) The quotient semigroup $J_a/S_{a)}$ is a completely (0)-simple semigroup.
3. $H_a \subset J_a$ is a regular $\mathcal{H}$-class if and only if $H_a$ is a maximal subgroup of $S$.

## 2.2 Review of the Rees Matrix Theorem

We provide the definition of Rees matrix semigroup and the statement of the Rees Matrix Theorem. Let $G$ be a group. Let $Q$ be an $m \times n$ matrix with entries in $G^0$. We define a *Rees matrix semigroup*, which we denote by $\mathcal{M}(G,Q)$ Let $\Lambda = \{1, \ldots, \lambda, \ldots, m\}$ index the rows of $Q$ and let $I = \{1, \ldots, i, \ldots, n\}$ index the columns of $Q$. An $n \times m$ monomial matrix $M$ with non-zero entry $g \in G$ in the $(r, \lambda)$ place will be denoted by $[r, g, \lambda]$. The $m \times n$ matrix with 0 in every place will be denoted by 0. The elements of the Rees matrix semigroup $\mathcal{M}(G,Q)$ are the monomial matrices over $G$, along with the 0 matrix. For $M, N$, monomial matrices over $G$, we define an associative multiplication: $M * N = MPN$; the 0 matrix acts as a 0. The reader can verify that the multiplication rule above is equivalent to the following multiplication given in terms of the 3-tuple form for monomial matrices.

- For $[r, g, \alpha]$ and $[s, h, \beta]$ ($r, s \in I$, $g, h \in G$ and $\alpha, \beta \in \Lambda$),
    - If $P(\alpha, s) \in G$, then $[r, g, \lambda][s, h, \mu] = [r, gP(\lambda, s)h, \mu]$.
    - if $P(\lambda, s) = 0$, then $[r, g, \lambda][s, h, \mu] = 0$.

**Example 2.5.** *Let $(X; G)$ be a group representation and let $P$ be a regular c-ramified matrix over $(X; G)$. Consider the semigroups of transformations $U((X; G); P)$ and $L((X; G); P)$. We use $V$ and $V_{\theta_P}$ to denote the ideals of constant functions of $U((X; G); P)$ and $L((X; G); P)$, respectively. It is not difficult to see that $U((X; G); P)/V$ and $L((X; G); P)/V_{\theta_P}$ are isomorphic to the same Rees matrix semigroup. In fact if we will define a matrix $Q$ with entries in $G^0$ in two steps. First if $P(\lambda, i) = x \in X$, let $Q(\lambda, i) = 0$ ($\lambda \in \Lambda$, $i \in I$). Second, regard the permutation entries of $P$ as elements from an abstract group. Abusing notation, let $G$ also denote the abstract group. Then $U((X; G); P)/V$ and $L((X; G); P)/V_{\theta_P}$ are both isomorphic to $\mathcal{M}(G; Q)$ with respective isomorphisms given by $[r, g, \lambda]/V \to [r, g, \lambda]$, $V \to 0$, and $[r, g, \lambda]/V_{\theta_P} \to [r, g, \lambda]$, where $V/\theta_P \to 0$.*

For a Rees matrix semigroup $\mathcal{M}(G,Q)$, the matrix $P$ is said to be a *sandwich matrix over $G$*. The sandwich matrix $Q$ is said to be *regular* if every row of $Q$ has at least one non-zero entry and every column of $Q$ has at least one non-zero entry. It is straight-forward to verify that $\mathcal{M}(G,Q)$ is completely (0)-simple if and only if $Q$ is regular, $\mathcal{M}(G,Q)$ is simple (in the semigroup sense) if and only if $Q$ has no 0 entries. For this reason, when we refer to a a Rees matrix semigroup, we assume that its sandwich matrix is regular. We let $J$ denote $\mathcal{M}(G,Q) \setminus \{0\}$. Note that in Example 2.5, the regularity of $P$ implies the regularity of $Q$.

**Remarks 2.6.** *Let $S = \mathcal{M}(G;Q)$ be a Rees matrix semigroup. It is easily verified that if $Q$ is an $m \times n$ matrix, then $J$ has $m$ distinct $\mathcal{R}$ classes and $n$ distinct $\mathcal{L}$ classes. In fact for $[r,g,\alpha]$ and $[s,h,\beta]$, we have $[r,g,\alpha] \mathrel{\mathcal{R}} [s,h,\beta]$ if and only if $r = s$ and $[r,g,\alpha] \mathrel{\mathcal{L}} [s,h,\beta]$ if and only if $\alpha = \beta$.*

*Thus we can use $I$ to simultaneously index the columns of $Q$ and the $\mathcal{R}$ classes of $J$ and use $\Lambda$ to simultaneously index the rows of $Q$ and the $\mathcal{L}$ classes of $J$. Moreover from Observation 2.1 and a result which comes later, Lemma 3.7, we have that if $(Y;J;\gamma)$ is a faithful representation of the maximal $\mathcal{J}$ class $J$ of $\mathcal{M}(G;Q)$, then $I$ also indexes the image sets of $\gamma(J)$ and $\Lambda$ also indexes the kernels of the functions of $\gamma(J)$.*

**Theorem 2.7.** *Let $S$ be completely (0)-simple. Then $S$ is isomorphic to a regular Rees matrix semigroup $\mathcal{M}(G,Q)$.*

We give a slightly non-standard proof of the finite Rees Matrix Theorem. For the purpose of the proof, we introduce some basic definitions. The definitions and notation that follow, particularly those involving what we call **Rees generating set** will be used throughout the paper.

**Definition 2.8.** *Let $S$ be completely (0)-simple semigroup with maximal $\mathcal{J}$ class $J$. Let $e$ be an idempotent of $J$. Let $H = eSe \cap J$. A subset $\mathcal{P} \subset J$ will be called a **Rees generating set based at e** if $\mathcal{P} = \mathcal{I} \cup \mathcal{K} \cup H$ where*

1. *$H = eSe \cap J$*

2. *$\mathcal{I}$ is a subset of $\mathcal{L}_e$ formed by choosing one representative from each $\mathcal{H}$ class of $\mathcal{L}_e$. We let $\mathcal{I} = \{f_r : r \in I\}$, where $I$ is an index set for the set of $\mathcal{R}$ classes of $J$ ( $I$ also indexes the $\mathcal{H}$ classes of $\mathcal{L}_e$). Observe that for all $r \in I$ we have that $f_r e = f_r$.*

3. *$\mathcal{K}$ is defined dually. That is, $\mathcal{K}$ is a subset of $\mathcal{R}_e$ containing exactly one representative from each $\mathcal{H}$ class of $\mathcal{R}_e$. We let $\mathcal{K} = \{g_\lambda : \lambda \in \Lambda\}$ where $\Lambda$ indexes the $\mathcal{L}$ classes of $J$. For all $\lambda \in \Lambda$, we have that $eg_\lambda = g_\lambda$.*

4. *If $\mathcal{P}$ is a Rees generating set and both $I$ and $\Lambda$ are totally ordered, then we say that $\mathcal{P}$ is an **ordered Rees generating set** (based at e).*

What does a Rees generating set $\mathcal{P}$ actually generate? The next lemma is a unique factorization theorem for elements of $J$ in terms of a given Rees generating set $\mathcal{P}$.

**Lemma 2.9.** *Suppose $S = J^0$ is a completely (0)-simple semigroup and that $e$ is an idempotent in $J$. Let $H = eSe \cap J$. Suppose $\mathcal{P} = \mathcal{I} \cup \mathcal{K} \cup H$ is a Rees generating set based at e. Then every element in $j \in J$ is uniquely represented as a product $j = fhg$ where $f$ is in $\mathcal{I}$, $h$ is in $H$ and $g$ is in $\mathcal{K}$.*

*Proof.* Note that $H$ is an $\mathcal{H}$ class of $J$. Let $\mathcal{P}$ be a Rees generating set based at $e$. We will claim that $J = \mathcal{I}H\mathcal{K} \setminus \{0\}$. Let $j \in J$. By definition of Rees generating set, there exists $r \in I$ and $\alpha \in \Lambda$ such that $j \mathcal{R} f_r$ and $j \mathcal{L} g_\alpha$ and $g_\alpha \in \mathcal{K}$, $f_r \in \mathcal{I}$. Thus there exists $v \in J$ such that $j = vg_\alpha$. Also, by Green's Lemma, $v \mathcal{R} j$. So $v \mathcal{R} f_r$. In particular, there exists $u \in J$ such that $v = f_r u$. We have $j = f_r u g_\alpha = (f_r e)u(eg_\alpha) = f_r(eue)g_\alpha$. Let $h = eue$ and observe that $j \in J$ implies that $eue \in J$. Thus, $eue \in H$. This completes the proof of the claim. Suppose also that $h' \in H$ such that $j = f_s h' g_\beta$, $f_s \in \mathcal{I}$, $h' \in G$ and $g_\beta \in \mathcal{K}$. Then by Green's Lemma, $f_s \mathcal{R} f_r$ and $g_\beta \mathcal{K} g_\alpha$. By definition of Rees generating set, we have $s = r$ and $\alpha = \beta$. By Corollary 2.3, item 2, we have that $h = h'$. □

*Proof.* (**Proof of finite Rees Matrix Theorem**) As remarked, if $Q$ is regular, then $\mathcal{M}(G; Q)$ is completely (0)-simple. We prove the converse direction.

Let $S$ be completely (0)-simple with maximal $\mathcal{J}$ class $J$. Choose an idempotent $e \in J$ and choose a Rees generating set $\mathcal{P}$ based at some idempotent $e \in J$. By Lemma 2.9, there exists a unique $\mathcal{P}$ decomposition of each element of $J$. This decomposition determines a bijection between elements of $J$ and triples of the form $\{[r, h, \alpha] : r \in I, h \in G, \alpha \in \Lambda\}$, the bijection given by $j = f_r h g_\alpha \mapsto [r, h, \alpha]$.

We define $P$, an $m \times n$ matrix with entries in $G^0$ as follows. For $r \in I$ and $\alpha \in \Lambda$, we let the $P(\alpha, r) = g_\alpha f_r$. Note that $P(\alpha, r) \in eJe$. By Corollary 2.3, item 1, we know that $H = eJe \cap J$ is a subgroup of $S$. So $P$ is as an $m \times n$ matrix with entries in $G^0$. We show that $P$ is regular. Suppose that $\alpha \in \Lambda$. We claim there exists $r \in I$ such that $g_\alpha f_r \neq 0$. By the definition of Rees generating set, we have that $g_\alpha e = g_\alpha$. But $g_\alpha \in J$ implies $g_\lambda e \in J$. By Green's Lemma, left multiplication by $g_\alpha$ maps $R_e$ bijectively to a subset of $J$. Since $R_e$ contains a member of $\mathcal{I}$, say $f_i$ ($I \in I$), we have that $g_\alpha f_i \neq 0$. Thus $P(\alpha, i) \neq 0$. By a similar argument, same we can show that for any $i \in I$ there exists $\beta \in \Lambda$ such that $P(\beta, i) \neq 0$. Thus $P$ is regular.

The reader can verify that the map $\kappa$ which sends elements of $J$ to their coordinates and which sends 0 to 0 is a semigroup isomorphism between $S$ and $\mathcal{M}(G; P)$. That $\kappa$ is a semigroup homomorphism follows immediately from the definition of the multiplication in a Rees matrix semigroup; $\kappa$ is one-to-one follows directly from the uniqueness of the factorization in Lemma 2.9. Thus, $S$ is isomorphic to $\mathcal{M}(G; P)$, a regular Rees matrix semigroup. □

We show that Rees generating sets can be used to describe any isomorphism from a completely (0)-simple semigroup into a Rees matrix semigroup. First we need the following useful lemma.

**Lemma 2.10.** *Let $\mathcal{M}(G; Q) = S$ be a regular Rees matrix semigroup. Then there exists a subgroup of $\mathcal{M}(G; Q)$ which is isomophic to $G$.*

*Proof.* Since $Q$ is not the 0-matrix, we can find $i \in I$ and $\lambda \in \Lambda$ such that $q_{\lambda, i} = g \in G$. Observe that $[r, g^{-1}, \lambda] = e$ is idempotent. By the second sentence of item 1 of Green's Lemma, $eJe \cap J = \{[i, h, \lambda] : h \in G\}$ is a subgroup of $S$. The map $[i, h, \lambda] \to gh$ is an isomorphism from $eJe \cap J$ to $G$. □

**Lemma 2.11.** *Let $S = \mathcal{M}(G; Q)$ be a finite regular Rees matrix semigroup. Then there exists $\mathcal{Q}$ an ordered Rees generating set for $J = S \setminus \{0\}$ such that $Q$ is exactly the matrix associated with $\mathcal{Q}$.*

*Proof.* As in the previous lemma, we choose a non-0 entry of $Q$ to locate a non-0 idempotent $e$ and consider the subgroup $eJe \cap J$. For $s \in I$, let $f_s = [s, Q(i, \lambda)^{-1}, \mu]$ and let for $\lambda \in \Lambda$, let $g_\lambda = [i, 1, \lambda]$. Letting $\mathcal{I}$ be the ordered set $\{f_r : r \in I\}$ and let $\mathcal{K}$ be the ordered set $\{g_\lambda : \lambda \in \Lambda\}$. Observe that $\mathcal{Q} = \mathcal{R} \cup \mathcal{K} \cup H$ is an ordered Rees generating set with associated matrix $Q$ (with respect to the isomorphism given in Lemma 2.10 from $eSe \cap J \to G$). □

## 3 R-representations

**Definition 3.1.** *Let $S$ be completely (0)-simple with maximal $\mathcal{J}$ class $J$. A representation $\gamma$ of $J$, $(Y; J; \gamma)$ will be called an **R-representation**.*

*An R-representation $(Y; J; \gamma)$ with the property that for all $j, k \in J$ if $jk = 0$, then $\gamma(j) \gamma(k)$ is a constant function, will be called a **c-representation**.*

Observe that any representation of a finite simple semigroup $S$ is a c-representation (since $S$ has no 0). Let $S$ be a semigroup with regular $\mathcal{J}$ class $J$. By Lemma 2.4, item 2c, if $a \in J$, then the regularity of $J = J_a$ implies that the quotient semigroup $J/S_{a)}$ is completely (0)-simple. Thus R-representations and the representations of the regular $\mathcal{J}$ classes of a finite semigroup are one and the same. To prove our primitivity results, we will only need to characterize a subclass of the c-representations defined below, the so-called **range-covered** c-representations. We provide a more general treatment of R-representations in Section 6.

### 3.1 Some additional definitions involving representations

The next series of definitions involves various weakenings of the transitivity property (a representation $(Y; U)$ is *transitive* if for all $a, b \in Y$ there exists $u \in U$ such that $u(a) = b$). Let $(Y; U)$ be a representation. We say that $(Y; U)$ is **range-covered** if $U(Y) = Y$. A transitive representation is obviously a range-covered representation. A representation $(Y; U)$ is said to be *cyclic* if there exists $y \in Y$ such that $Y = U(\{y\}) \cup \{y\}$. Cyclic representations are of interest in automata theory. Transitive and cyclic representations of finite simple semigroups were characterized by ([Stoll], [Tully]) in some of the earliest papers of algebraic semigroup theory.

The next series of definitions deal with pairs of elements of the carrier set which behave identically under all the operations. Let $(Y; U)$ be a representation. An equivalence relation $\alpha$ on $Y$ is a **deflationary** relation if for all $a, b \in Y$ and all $u \in U$, $u(a) = u(b)$.

Observe that if $\alpha$ is a deflationary relation, then $(Y; U)$ defines a representation of $U$ on $Y/\alpha$ which we denote by $(Y_\alpha; U)$. We say that $(Y_\alpha; U)$ is a **deflation** of $(Y; U)$. Note that the diagonal relation (which we denote by $\Delta$) is a deflationary relation and that $(Y; S)$ has a unique maximal deflationary relation which we call $\delta$. With $\delta$ the maximal deflationary equivalence, $(Y_\delta; S)$ has as its maximal deflationary relation the diagonal relation (on $Y/\delta$). If a representation $(Y; U)$ has as its maximal deflationery relation the diagonal relation, we say that $(Y; U)$ is **reduced**. Note that if $|Y| > 2$, then every primitive representation $(Y; S)$ is a reduced representation.

We borrow the next piece of terminology from tame congruence theory. For a representation $(Y; S)$, we say that a subset $U$ of $Y$ is a *neighborhood* of $Y$ if $U = e(Y)$ where $e \in S$ and $e^2 = e$. A neighborhood $U$ of $Y$ is a *minimal neighborhood* if $|U| > 1$ and if $V$ is a neighborhood properly contained in $U$, then $|V| = 1$. One consequence of Theorem 1.8 is that if $(Y; S)$ is primitive, then the minimal sets and the minimal neighborhoods of $(Y; S)$ are the same.

## 3.2 Observations concerning R-representations

I

**Example 3.2.** *Let $P$ be a c-ramified $m \times n$ matrix over the group representation $(X; G)$. Consider the action $(V_{\theta_P}; J((X; G); P)$. As remarked in Example 2.5, the semigroup $L((X; G); P)/V_{\theta_P}$ is a Rees matrix semigroup, so by the Rees Matrix Theorem, it is a completely 0-simple (with maximal $\mathcal{J}$ class $J((X; G); P)$). Thus $(V_{\theta_P}; J((X; G); P)$ is a c-representation.*

Let $(Y; S; \gamma) = (Y; S)$ be any finite representation and let $e$ be an idempotent of $S$. Since $J_e$ is a regular $\mathcal{J}$ class of $S$, we have that $eJ_e e \cap J_e$ is a subgroup and the representation $(e(Y); eJ_e e \cap J_e)$ is a group representation. Hence $(Y; S)$ is an extension of the group representation $(e(Y); eJ_e e \cap J_e)$. The next lemma is predicted by [3].

**Notation 3.3.** *Given an R-representation $(Y; J; \gamma)$, if we suppress $\gamma$, as we have mentioned for $j \in J$, $y \in Y$, we abbreviate $\alpha(j)(y)$ to $j(y)$. Note that if $j, k \in J$ and $jk \neq 0$, then $(jk)(y) = j(k(y))$ (if $jk = 0$, then $(jk)(y)$ is not even defined). More generally, if $j_1, \ldots, j_k \in J$ and $j_1 \ldots j_k \neq 0$, then it is acceptable to re-parenthesize $(j_1 \ldots j_k)(y)$ in many ways. In such situations, we re-parenthesize as is convenient without comment, including eliminating or adding parentheses.*

**Lemma 3.4.** *(Minimal Neighborhoods Lemma) Let $(Y; J)$ be an R-representation. Let $e, f \in J$ be idempotents. Then $(e(Y); eJe \cap J)$ is equivalent to $(f(Y); fJf \cap J)$.*

*Proof.* Let $e, f$ be idempotents of $J$. We have that $eJe \cap J$ and $fJf \cap J$ are subgroups; hence, both $(e(Y); eJe \cap J)$ and $(f(Y); fJf \cap J)$ are group representations. Because $e$ and $f$ are $\mathcal{J}$ related and $J$ is regular, there exist $s, t \in J$ such that $set = f$. Since $f$ is idempotent, $setset = f = set$. From Lemma 2.2, $etse$ is $\mathcal{H}$ related to $e$. By Lemma 2.3, item 2, $setset = set$ implies $etse = e$.

Let $\beta$ be the map which sends $e(y) \to s(e(y))$ $(e(y) \in e(Y))$; let $\beta$ send $eje \to sejet$ $(eje \in eJe \cap J)$. Note that if $eje, eke \in eSe \cap J$ and $(eje)(eje) \neq 0$, then $\beta(eje)\beta(eke) = (sejet)(seket) = sej(etse)(ket) = sej(e)ket = s(ejeeke)t$
$= \beta((eje)(eke))$.

Note that $set(Y) \subset se(Y)$; that is, $f(Y) \subset se(Y)$. But $se = setse$; hence, $se(Y) = f(Y)$. By Lemma 2.3, item 2, $s(eJe \cap J)t = fJf \cap J$. That is $\beta$ maps $e(Y) \cup eJe \cap J$ into $f(Y) \cup fJf \cap J$. By reversing the roles of $e$ and $f$, we can construct a map $\beta'$ which carries $f(Y) \cup (fJf \cap J)$ onto $e(Y) \cup (eJe \cap J)$. By finiteness, $\beta$ is a bijection.

Let $e(y) \in e(Y)$ and let $eje \in eJe \cap J$. Then $\beta(eje)(e(y)) = s(eje)(e(y)) = sejetse(e(y)) = sejet(s(e(y))) = \beta(eje)(\beta(e(y)))$. We have proven that $\beta$ is an equivalence. □

**Lemma 3.5.** *Let $(Y; J)$ be a c-representation. For the representation $(Y; J^Y)$, for all $j \in J$ we have that $j(Y)$ is a minimal neighorhood.*

*Proof.* Since $J$ is a finite regular $\mathcal{J}$ class, for every element $j \in J$ there exists an idempotent $e \in J$ such that $e \mathcal{R} J$. Hence by Observation 2.1, $e(Y) = j(Y)$. Thus every image set is an idempotent image. If $(Y; J)$ is a c-representation, then by definition of c-representation the composition of two representations of elements of $J$ is either represented by an element of $J$ or the composition results in a constant function. Thus among the non-constant images of the representations of $J^Y$, the set $j(Y)$ is minimal. □

Let $(Y; U; \gamma)$ be a representation, where $U$ is a subset of a semigroup $S$. Suppose that $Y \subset Z$ and $U \subset V \subset S$. If we can extend $\gamma$ to $\tilde{\gamma}$ so that $(Z; V; \tilde{\gamma})$ is a representation and for all $u \in U$, $\tilde{\gamma}(u)|_Y = \gamma(u)$ then we say that $(Z; V; \tilde{\gamma})$ is an **expansion** of $(Y; U; \gamma)$. If $W \subseteq Y$ and $\gamma(U)$ maps $W$ into itself, then we refer to $(W; U; \gamma|_W)$ as a **subrepresentation** of $(Y; U; \gamma)$.

**Definition 3.6.** *If $(Y; S)$ is a representation and $X$ is a minimal neighborhood of $(Y; S)$ with $e^2 = e \in S$ be such that $e(Y) = X$, then we say that $(Y; S)$ is an **m-expansion** of the representation $(X; eJe \cap J)$.*

By the Minimal Neighborhoods Lemma, if $(Y; J; \gamma) = (Y; J)$ is a c-representation, then it is the m-expansion of a unique up to equivalence group representation $(e(Y); eJe \cap J)$ where $e$ is any idempotent of $J$.

The next lemma provides necessary and sufficient conditions for faithfulness of an R-representation $(Y; J; \gamma)$. We will need to be explicit about $\gamma$. Let $e$ be an idempotent of a $\mathcal{J}$ class $J$. We let $(e(Y); eJe \cap J; \gamma|_{e(Y)})$ be the representation with $\gamma|_{e(Y)} : eJe \cap J \to S_{e(Y)}$, $eje \to \gamma(eje)|_{e(Y)}$ ($eje \in eJe \cap J$).

**Lemma 3.7.** *An R-representation $(Y; J; \gamma)$ is faithful if and only if*

1. *$(e(Y); eJe \cap J; \gamma|_{e(Y)})$ is faithful.*

2. *Two elements of $j, k \in J$ are $\mathcal{R}$ related if and only if $im(\gamma(j)) = im(\gamma(k))$.*

3. *Two elements are $j, k \in J$ are $\mathcal{L}$ related if and only if $ker(\gamma(j)) = ker(\gamma(k))$.*

*Item 2 holds if and only if the number of minimal neighborhoods is the same as the number of $\mathcal{R}$ classes of $J$.*

*Proof.* In this proof we will need to distinguish between elements of $J$ and their representation under $\gamma$. Suppose that $(Y; J; \gamma)$ is faithful. Let $j, k \in J$ be such that $im(\gamma(j)) = im(\gamma(k))$. We showed in Lemma 3.5 that there exists an idempotent $e \in J$ such that $im(\gamma(e)) = im(\gamma(j))$. We have $\gamma(j) = \gamma(ej) = \gamma(e)\gamma(j)$ and $\gamma(k) = \gamma(ek) = \gamma(e)\gamma(k)$. Since $\gamma$ is faithful, $ej = j$ and $ek = k$; hence, by the second sentence of item 1 of Green's Lemma, it follows that $j \mathcal{R} k$. Also $ker(j) = ker(k)$ implies that $j \mathcal{L} k$ can be similarly be argued.

Suppose that $e \in J$ is idempotent and for some $j, k \in J$ and $eje \neq eke$. We show that $\gamma|_{e(Y)}(eje) \neq \gamma|_{e(Y)}(eke)$. Because $(Y; J; \gamma)$ is faithful, $\gamma(eje) \neq \gamma(eke)$. Thus there exists $y \in Y$ such that $\gamma(eje)(y) \neq \gamma(eke)(y)$; hence, $\gamma(eje)(e(y)) \neq \gamma(eke)(e(y))$. It follows that $\gamma|_{e(Y)}(eje) \neq \gamma|_{e(Y)}(eke)$. We have shown that if items 1-3 hold, then $(Y; J; \gamma)$ is faithful.

For the converse assume that items 1-3 hold for the R-representation $(Y; J; \gamma)$. Let $j, k \in J$ be such that $\gamma(j) = \gamma(k)$. By items 2 and 3, we have that $j \mathcal{H} k$. Let $e \in J$ be an idempotent. Since $j \mathcal{J} e$, there exist $s, t \in J$ such that $sjt = e$. By item 2 of Green's Lemma, we have $skt \mathcal{H} sjt = e$. Note also that $esjte = sjt$ and $eskte = skt$ and that $\gamma(j) = \gamma(k)$ implies that $\gamma(esjte) = \gamma(eskte)$. So $esjte, eskte \in eJe \cap J$ and from the previous sentence, we have $\gamma|_{e(Y)}(eskte) = \gamma|_{e(Y)}(esjte)$. From item 1, we have $eskte = esjte$. Thus $skt = sjt$. Since $j \mathcal{H} k$, by Lemma 2.3, item 2, we have $j = k$. We have shown that items 1-3 imply that $\gamma$ is faithful.

By Observation 2.1, the number of $\mathcal{R}$ classes of $\mathcal{J}$ is an upper bound for the number of minimal neighborhoods of $(Y; J; \gamma)$. Thus if the number of minimal neighborhoods is equal to the number of $\mathcal{R}$ classes, then item 2 obviously holds, and conversely. $\square$

## 3.3 Ramified Rees representations

We define the **ramified** matrices over a group representation $(X; G)$.

**Definition 3.8.**    1. Let $(X; G)$ be a group representation and let $N_X = T_X \setminus S_X$ (i.e., $N_X$ is the set of non-permutation transformations of $X$). An $m \times n$ matrix $P$ with entries in $G \cup N_X$ will be said to be a **ramified matrix** over $(X; G)$.

2. Let $(V; J((X; G); P)$ denote the action $v \to MPv$ ($v \in V$, $M \in J((X; G); P)$).

3. Let $G$ be an abstract group, $Q$ be a $m \times n$ matrix over $G^0$, and $(X; G; \gamma)$ be a representation of $G$. Then a ramified matrix $P$ over $(X; G; \gamma)$ will be said to be a **ramification of Q** over $(X; G; \gamma)$ if $Q(\lambda, i) = g \in G$ implies that $P(\lambda, i) = \gamma(g)$ and if $Q(\lambda, i) = 0$, then $P(\lambda, i) \in N_X$.

If all the non-permutation entries of $P$ are constant functions (i.e. in $X$), then $P$ is a c-ramified matrix over $(X; G)$ (see Definition 1.1). The only change from Definition 1.2 is that $P$ is only required to be a ramified matrix, rather than c-ramified as in Definition 1.2.

We widen the set of examples of R-representations associated with a given ramified Rees matrix to what we call the **ramified Rees representations**. Every R-representation will be turn out to be a variant of a ramified Rees representation. Let $(V; J((X;G); P))$ be a representation as defined in Definition 3.8. Let $\alpha$ be a deflation of $V$. We denote $V/\alpha$ by $V_\alpha$.

**Definition 3.9.** Let $(X; G; \gamma)$ be a group representation and let $P$ be a regular ramified matrix over $(X; G)$.

- A **ramified Rees representation** is any representation of the form $(V_\alpha; J((X;G); P)$ where $\alpha$ is a deflationary relation.

- If the entries of $P$ are contained in $X \cup \gamma(G)$ then we say that $(V_\alpha; J((X;G); P)$ is a **c-ramified Rees representation**.

Let $\delta$ denote the maximal deflation of $(V; J((X;G); P))$. Not surprisingly, $\delta$ can be described directly from $P$.

**Lemma 3.10.** Let $r_x$ and $s_y$ be contained in $V$.

1. $(r_x, s_y) \in \theta_P$ if and only if for all $(i, g, \lambda) \in J((X;G); P)$, $(i, g, \lambda)(r_x) = (r, g, \lambda)(s_y)$. Thus $\theta_P = \delta$, the maximal deflation of $(V; J((X;G); P))$.

2. For any $i \in I$ and any $x, y \in X$, we have that $(i_x, i_y) \in \theta_P$ implies that $x = y$.

*Proof.* In Definition 1.4, we have remarked that $(r_x, s_y) \in \theta_P$ if and only if for all $\lambda \in \Lambda$, $P(\lambda, r)(x) = P(\lambda, s)(y)$ $(x, y \in X, r, s \in I)$. Item 1 follows directly from the definition of the action of $J((X;G); P)$ on $V$.

For item 2 suppose that $(i_x, i_y) \in \theta_P$. Then for all $\lambda \in \Lambda$ we have $P(\lambda, i)(x) = P(\lambda, i)(y)$. Since $P$ is regular, there exists $\lambda \in \Lambda$ such that $P(\lambda, i)$ is a permutation of $X$. It follows that $x = y$. □

### 3.3.1 Fundamental lemma of R-representations

We say that $(Y; J; \gamma)$ is **m-faithful** if it m-expands a faithful group representation.

**Lemma 3.11.** *(Fundamental Lemma of R-representations)* Let $J^0$ be a completely 0-simple semigroup with R-representation $(Y; J)$.

1. If $(Y; J)$ is a m-faithful, range-covered R-representation, then $(Y; J)$ is equivalent to a ramified Rees representation.

2. If $(Y; J)$ is a m-faithful range-covered c-representation, then $(Y; J)$ is equivalent to a c-ramified Rees representation.

3. Furthermore, if $J^0 = \mathcal{M}(G; Q)$, a Rees matrix semigroup, then in item 1 and item 2 above, we have that $(Y; J)$ is equivalent to $(V_\alpha; J((X;G); P))$ where $P$ is a ramification of $Q$ over $(X; G)$ and $(X; G)$ is equivalent to *the unique up to equivalence group representation for which $(Y; J)$ is an m-expansion.*

*In the case that $(Y; J)$ is a c-representation, $P$ is a c-ramification of $Q$ over $(X; G)$.*

*Proof.* (Proof of Fundamental Lemma of R-representations)

Let $(Y; J; \gamma)$ be an m-faithful R-representation which is range covered with idempotent $e \in J$. By Lemma 2.9, we can choose $\mathcal{Q} = \mathcal{I} \cup \mathcal{K} \cup H$, a Rees generating set based at an idempotent $e \in J$. In fact, by the Rees Matrix Theorem, $J^0$ is isomorphic to a Rees matrix semigroup $\mathcal{M}(G; Q)$ under a map which we will call $\alpha$. We assume that $\mathcal{Q}$ is the image under $\alpha$ of a Rees generating set (of the maximal $\mathcal{J}$ class of $\mathcal{M}(G; Q)$) which returns the matrix $Q$. Such a Rees generating set exists by Lemma 2.11.

To keep the notation simple, for $j \in J$, $y \in Y$, we refer to $\gamma(j)(y)$ as $j(y)$. Let $H = eJe \cap J$, a subgroup of $J^0$ which is isomorphic to $G$ by Lemma 2.10. Let $e(Y) = X$. By construction we have that $H$ maps $X$ to itself. Hence, $H$ can be represented on $X$ via $(X; H)$ $(=(X; H; \gamma|_X))$. Since $(Y; J)$ is m-faithful, from the Minimal Neighborhoods Lemma, it follows that that $(X; H)$ is faithful.

For $i \in I$ and $\lambda \in \Lambda$, observe that $g_\lambda f_i$ is contained in $eJe$. In particular, $g_\lambda f_i$ maps $X$ to itself; hence, $g_\lambda f_i$ induces a function from $X$ to itself.

We denote this function on $X$ by $(g_\lambda f_i)|_X$. Recall that $N_X$ denotes $T_X \setminus H|_X$. We have that $(g_\lambda f_i)|_X$ is an element of $N_X \cup H|_X$. We define $P$, an $m \times n$ ramified matrix over $(X; H)$. For $\lambda \in \Lambda$ and $i \in I$ let $P(\lambda, i) = (g_\lambda f_i)|_X$. Note that $P(\lambda, i)$ is contained in $N_X \cup H|_X$. Since $Q$ is regular, it follows that $P$ is regular.

We define a function $\kappa : J \mapsto J((X; G); P)$, where $J((X; G); P)$ is the set of monomial $n \times m$ matrices over $(X; H; \gamma|_X)$ as follows. By Lemma 2.9, each $j \in J$ has a unique $\mathcal{Q}$ factorization $j = f_r h g_\lambda$ ($r \in I$, $h \in H$, $\lambda \in \Lambda$). Let $\kappa(j) = [r, h|_X, \lambda]$, the $m \times n$ monomial matrix with a non-zero entry $h|_X$ in the $(r, \lambda)$ entry. Since $(X; H)$ is faithful, it follows that $\kappa$ is one-to-one on $J$. Moreover, for $r \in I$, $\lambda \in \Lambda$ and $h \in H$, $[r, h|_X, \lambda] \in J((X; G); P)$ is the image of $f_r h g_\lambda$. Hence, $\kappa(J) = J((X; G)); P)$. We extend $\kappa$ to $Y$ so that it maps $Y$ to $V(m; X) = V$.

We define $\kappa$ on $X = \gamma(e)(Y)$ first. For $x \in X$ let $\kappa(x)$ be defined as $1_x$. Next let $y \in Y$ be arbitrary. By the range-covered hypothesis, $J(Y) = Y$; hence, there exists $j \in J$ such that $y \in j(Y)$. By definition of Rees generating set, there exists a unique $f_r \in \mathcal{I}$ such that $f_r$ and $j$ are in the same $\mathcal{R}$ class. By Observation 2.1, $f_r(Y) = j(Y)$. By definition of Rees generating set based at $e$, $f_r e = f_r$; hence, $f_r e(Y) = j(Y)$. Since $f_r e \in J$ and $\gamma$ is an R-representation, for all $z \in Z$, we have $(f_r e)(z) = f_r(e(z))$. Thus there exists $x \in X$ such $y = f_r(x)$. We define $\kappa(y)$ to be $r_x$.

At this point the restriction of $\kappa$ to $Y$ may be multi-valued. We turn it into a single-valued function by making some identifications. For $x, w \in X$ and $r, s \in I$ let $(r_x, s_w) \in \alpha$ if $f_r(x) = f_s(w)$. If $y = f_r(x) = f_s(w)$, then for all $\lambda \in \Lambda$, we have $g_\lambda f_r(x)) = g_\lambda f_s(w))$. Hence if $(r_x, s_w) \in \alpha$, then for all $\lambda \in \Lambda$, $P(\lambda, r)(x) = P(\lambda, s)(w)$. In particular, $(r_x, s_w) \in \theta_P$. Thus, $\alpha \leq \theta_P$. We define a mapping from $(Y;J)$ into $(V_\alpha; J((X;G);P))$. We denote this map by $\kappa_\alpha$, where $\kappa_\alpha : (Y;J) \to (V_\alpha; J((X;G);P))$ and $\kappa_\alpha(y) = \kappa(y)/\alpha$ ($y \in Y$) and $\kappa_\alpha(j)$ $= \kappa(j)$ ($j \in J$).

By the definition on $\kappa_\alpha$, we have $\kappa_\alpha$ is one-to-one on $Y$. We have mentioned that $\kappa$ maps $J$ bijectively to $J((X;G);P)$, so $\kappa_\alpha$ maps $J$ bijectively to $J((X;G);P)$. Let $r_x/\alpha \in V_\alpha$. Since $r_x/\alpha$ has pre-image $\gamma(f_r)(x)$ under $\alpha$, $\kappa_\alpha|_Y$ is a bijection between $Y$ and $V_\alpha$. Thus $\kappa_\alpha$ maps $Y \cup J \to V_\alpha \cup J((X;G);P)$ bijectively.

We prove that $\kappa_\alpha : (Y;J) \to (V_\alpha; J((X;H);P))$ is an equivalence. As a first step, we show first that for all $y \in Y$ and $j \in J$, we have $\kappa(j(y)) = \kappa(j) \kappa(y)$. Let $j$ have unique $\mathcal{Q}$ decomposition $f_r h g_\lambda$. As we have shown, there exists $i \in I$ and $x \in X$ such that $y = f_i(x)$. For any such choice of $i$ and $x$, we have $j(y) = f_r h g_\lambda(f_i(x))$. Since $f_r h g_\lambda \neq 0$ and $\gamma$ is an R-representation, we have $f_r h g_\lambda(f_i(x)) = f_r) h g_\lambda(f_i(x))$. Thus, $\kappa(j(y)) = r_{h g_\lambda(f_i)(x)}$. If $g_\lambda f_i \neq 0$, then since $\gamma$ is an R-representation, we have $h g_\lambda(f_i)(x)) = h g_\lambda f_i(x) = h|_X P(\lambda, i)(x)$, the last identity following from definition of the matrix $P$ and the fact that $g$ maps $Y$ to $X$. On the other hand, if $g_\lambda f_i = 0$, then $g_\lambda f_i$ is a constant function with singleton range $\{x'\}$, where $x' \in X$. Note in this case also that $x' = P(\lambda, i)$. Thus, $h g_\lambda(f_i(x)) = h(P(\lambda, i))(x) = h|_X(P(\lambda, i))(x)$. For both cases $h g_\lambda(f_i)(x)) = h|_X(P(\lambda, i))(x)$. Thus for all $j \in J$, $y \in Y$, $\kappa(j)(y)) = r_{h|_X(P(\lambda,i)))}/\alpha$. On the other hand, by the definition of the action in $(V; J((X;H);P))$, we have $\kappa(j) (\kappa(y)) = [r, h|_X, \lambda](i_x) = r_{h|_X,P(\lambda,i)(x)}$. Note that since $\alpha$ is a deflation, if we had chosen a different set of coordinates for $y =$, say $s_y$, then $(i_x, s_y) \in \alpha$ implies that $\kappa(j)(i_x) = \kappa(j)(s_y)$. Thus for $y \in Y$, $j \in J$, we have shown that $\kappa(j(y)) = \kappa(j)$. From the definition of $\kappa_\alpha$ and the fact that $\alpha$ is a deflation, it now follows with a moment's thought that for $j \in J$ and $y \in Y$ to $J$, $\kappa_\alpha(j(y)) = \kappa_\alpha(j)\kappa_\alpha(y)$.

To show for any $j, k \in J$, $jk \neq 0$, $\kappa_\alpha(jk) = \kappa_\alpha(j)\kappa_\alpha(k)$, arguments that are similar in nature to our Rees generating set proof of the Rees Matrix Theorem can be used. We leave it to the reader to check the computations. This completes the proof of item 1. That is, $(Y;J;\gamma)$ is equivalent to a ramfied Rees matrix semigroup $(V_\alpha; J((X;H);P))$.

To prove item 2, it is easy to see that $(Y;J;\gamma)$ above is a c-representation if and only if the matrix $P$ is a c-ramified matrix. For item 3 of the theorem, at the beginning of the proof of item 1, we chose the Rees generating set $\mathcal{Q}$ in such a way that the resulting matrix $P$ is a ramification of $Q$ over $(X;H)$.

Since $H$ is isomorphic to $G$, it follows that in one more step we can say that $P$ is a ramification of $Q$ over $(X;G)$. □

**Observation 3.12.** *The Fundamental Lemma of R-representations provides a clear-cut way to produce, up to equivalence, every m-faithful range-covered representation of a completely 0-simple semigroup $S$. First find an isomorphism between $S$ and a Rees matrix semigroup $\mathcal{M}$ $(G;M)$ and a faithful representation $(X;G)$ of $G$. Consider the set of all ramifications of $M$ over $(X;G)$. Each ramification $P$ leads to a set of m-faithful range covered R-representations, that set ranging over the equivalence relations $\alpha$ contained in $\theta_P$.*

We provide necessary and sufficient conditions for faithfulness of an R-representation $(V_\alpha; J((X;G); P))$ in terms of the parameters $\alpha$, $P$, and $(X;G)$.

**Definition 3.13.** *For the representation $(V_\alpha; J((X;G); P))$, for $i \in I$, we let $X_i = \{i_x/\alpha : x \in X\}$. Observe that for all $i \in I$, by Lemma 3.10, item 2, we have $|X_i| = |X|$.*

**Lemma 3.14.** *Let $(X;G)$ be a group representation, let $P$ be a regular ramified $m \times n$ matrix over $(X;G)$ and let $\alpha \leq \theta_P$. Consider the ramified Rees representation $(V_\alpha; J((X;G); P))$. For any $[r, g, \mu], [s, h, \beta] \in J((X;G); P)$, $\mu \neq \beta$, the kernels of $[r, g, \mu]$ and $[s, h, \beta]$ are the same if and only if $P_\mu$ is proportional to $P_\beta$.*

*In particular, if $(V_\alpha; J((X;G); P))$ is faithful, then $P$ is right reductive.*

*Proof.* Assume that $(X;G)$ is faithful. If $P_\mu$ and $P_\beta$ are proportional ($\mu, \beta \in \Lambda$), then it is straightforward to verify that the kernels of $[a, h, \mu]$ and $[b, g, \beta]$ are the same ($a, b \in I$, $h, g \in G$).

Suppose $P_\mu$ and $P_\beta$ are not proportional. If there exists $i \in I$ such that $P(\mu, i) \in X$ but $P(\beta, i) \in G$, then then $[a, h, \mu]$ identifies the set $X_i = \{i_x/\alpha : x \in X\}$; whereas, the fact that $P(\beta, i) \in G$ it follows that $[b, h, \beta]$ maps $X_i$ bijectively. But by the Minimal Neighborhoods Lemma, we have $|X_i| = |X| > 1$. Thus for this case, $ker([a, h, \mu]) \neq ker([b, h, \beta])$. So we may as well assume that $P(\mu, i) \in X$ if and only if $P(\beta, i) \in X$ ($i \in I$). Since $P$ is regular, there exists $i \in I$ such that $P(\mu, i) \in G$. From the first paragraph of this proof, multiplying $P_\mu$ on the left by a group group element, multiplying $P_\beta$ on the left by a group element do not effect the kernels of $[a, g, \mu]$, $[b, h, \beta]$ respectively. Thus without loss of generality we may assume that $P(\mu, i) = P(\beta, i) = 1$, the identity function of $X$. That $P_\mu$, $P_\beta$ are not proportional now implies that there exist $r \in I$, $y \in X$ such that $P(\mu, r)(y) \neq P(\beta, r)(y)$. Let $P(\mu, r)y = x \in X$. We have $a_{gP(\mu,r)(y)} = [a, g, \mu](r_y) = [a, g, \mu](i_x) = a_{g(x)}$. But now $[b, h, \beta](i_x) = b_{h(x)}$. Whereas, $P(\beta, r) \neq x$ implies that $[b, h, \beta](r_y) = b_{hP(\beta,r)(y)} \neq b_{h(x)}$, the last non-equality holding since $h$ is a permutation of $X$. In particular, the kernels of $[a, g, \alpha]$ and $[b, h, \beta]$ are not the same. □

**Lemma 3.15.** *Let $(X;G)$ be a group representation, let $P$ be a regular ramified $m \times n$ matrix over $(X;G)$ and let $\alpha \leq \theta_P$. Consider the ramified Rees representation $(V_\alpha; J((X;G); P))$. We have $(V_\alpha; J((X;G); P))$ is faithful if and only if the following hold.*

1. *$(X;G)$ is faithful.*

2. *$P$ is right reductive.*

3. *The number of neighborhoods of $(V_\alpha; J((X;G); P))$ equals the the number of columns of $P$ $(= n)$.*

*Moreover, $(V_\alpha; J((X;G); P))$ is faithful and reduced if and only if $\alpha = \theta_P$ and $P$ is reductive.*

*Proof.* The necesssity and sufficiency of the three items above follows from Lemma 3.7 and Lemma 3.14. For the last paragraph of the lemma, observe that $[r, g, \mu]$ and $[s, h, \beta]$ have the same image if and only if $X_r = X_s$ if and only if there exists $k \in S_X$ such that for all $x \in X$, $(r_{k(x)}, s_x) \in \alpha$. In fact since $P$ is regular, it follows that $k \in G$. In particular, if $P$ is left reductive if and only if the number of minimal neighborhoods of $(V_{\theta_P}; J((X; G); P))$ is $m$. Thus, $(V_{\theta_P}; J((X; G); P))$ is faithful if and only if $P$ is reductive. Earlier in Lemma 3.10 we saw that $(V_\alpha; J((X; G); P))$ is reduced if and only if $\alpha = \theta_P$. □

## 4 Proof of the main theorem

### 4.1 Elementary facts concerning primitive representations

Let $Y$ be a set and let $\theta$ be a binary relation of $Y$. We let $Eq(\theta)$ be the equivalence relation generated by $\theta$. Now suppose that $(Y; S)$ is a representation and that $\theta$ is once again a binary relation of $Y$. We let $Cg(\theta)$ denote the S-compatible equivalence relation generated by $\theta$. If $\theta = \{(a, b)\}$ $(a, b \in Y)$, we let $Cg(a, b)$ denote $Cg(\theta)$.

If $W \subseteq Y$, then we let $Cg(W)$ denote the smallest S-compatible equivalence on $Y$ which identifies $\{(u, v) : u, v \in W\}$ (rather than the more proper $Cg(W^2)$). The following lemma is a version for representations of what as known as Mal'cev's Lemma.

**Lemma 4.1.** *Let $(Y; S)$ be a representation and let $U$ be a binary relation on $Y$. Then $Cg(U) = Eq(\{(s(c), s(d)) : s \in S^1, (c, d) \in U\})$.*

From Malcev's Lemma we have : For $a, b, c, d \in S$, $(a, b) \in Cg(c, d)$ if and only if there exists a sequence $a = a_1, \ldots, a_{n+1} = b$ such that for $1 \leq i \leq (n)$ there exists $s_i \in S^1$ such that $\{a_i, a_{i+1}\} = \{s_i(c), s_i(d)\}$.

**Lemma 4.2.** *Suppose that $(Y; S)$ is a representation.*

1. *If $|Y| = 2$, then $(Y; S)$ is primitive*

2. *$(Y; S)$ is primitive if and only if for every pair $c \neq d \in Y$, $Cg(c, d) = \nabla$.*

### 4.2 R-representations and primitive representations

In this subsection we will also assume that $(Y; S)$ is primitive and in view of Lemma 4.2, item 1, we will assume that $|Y| > 2$. Recall that a subset $U \subset Y$ is said to be a minimal set if among all non-singleton image sets of $S$, $U$ is minimal with respect to set inclusion. Let $M(Y; S)$ denote the minimal sets. We say that $s$ is a *minimal function* if the image of $s$ is a minimal set. As we have observed, if $|Y| > 2$, then $M(Y; S)$ is non-empty.

Departing from our usual convention of letting $J$ refer to a regular $\mathcal{J}$ class, we let $J$ denote the set of minimal functions of a representation $(Y; S)$. Under the $|Y| > 2$ assumption we prove that if $(Y; S)$ is a faithful finite primitive representation, then the set of minimal functions $J$ is indeed a regular $\mathcal{J}$ class. It then follows that $(Y; J)$ is a c-representation. Let $C$ denote the (possibly empty) set of constant functions of $S$.

Let $\rho = Eq(R)$ where R = $\{U^2 : U \in M(Y; S)\}$). If $|Y| > 2$, then since there exist minimal sets, observe that $\rho \neq \Delta$. Much of the following lemma is predicted by [3].

**Lemma 4.3.** *Let $(Y; S)$ be a finite primitive representation, $|Y| > 2$ and let $J$ denote the minimal functions of $(Y; S)$.*

1. *If $s, t \in S^1$ are both non-constant, then $sS^1t$ contains a non-constant function.*

2. *$J$ is a $\mathcal{J}$ class of $S$.*

3. *$J$ is regular.*

4. *In the poset of $\mathcal{J}$ classes of $S$, the $\mathcal{J}$ class $J$ is the unique cover of the $\mathcal{J}$ class $C$ of left-0's of $S$ ($C$ might be empty).*

5. *$\rho = \nabla$*

6. *$(Y; J)$ is a range-covered c-representation.*

*Proof.* We prove item 1. Suppose $s, t \in S$ are both non-constant. Verify that there exist $c, d \in Y$ such that $s(c) \neq s(d)$ and $t(c) \neq t(d)$. Let $s(c) = a$ and $s(d) = b$. Since $(Y; S)$ is primitive, we have $(c, d) \in Cg(a, b)$. By Malcev's Lemma, there exist $f_1 \ldots f_{n-1}$ ($n - 1 \geq 1$) contained in $S$, and a sequence of elements in $Y$, $c = a_1 \ldots a_n = d$, such that $\{f_i(a), f_i(b)\} = \{a_i, a_{i+1}\}$ ($i = 1 \ldots n - 1$). Apply $t$ to each element of the chain above. Since $t(c) \neq t(d)$, it follows that for at least one of $i = 1 \ldots n - 1$, $tf_is(c) \neq tf_is(d)$. Hence $tf_is$ is non-constant. This proves item 1.

We prove items 2 and 3. If $j \in J$, then by the paragraph above, $jS^1j$ contains a non-constant. Let $s \in S^1$ such that $jsj$ is not constant; $j(Y)$ is a minimal set and $jsj$ not constant together imply that $j(Y) = jsj(Y) = (js)(j(Y))$. Thus there exists a positive integer $n$ such that $(js)^n = e$ is idempotent and non-constant. Clearly $e(Y) = j(Y)$; hence $e \in J$. Also we have that $ej = j$. It follows from $ej = j$ and $(js)^n = e$ that $e$ and $j$ are in the same $\mathcal{J}$ class. Let $j, k$ be arbitrary element of $J$. By two applications of the previous paragraph there exist $s, t \in S^1$ such that $jsktj$ is not constant. Since $jsktj$ is a minimal function, it follows that there exists a positive integer $m$ such that $(jskt)^m = f$ is idempotent, minimal and $fj = j$. From $(jskt)^m = f$, we have $f \leq_J j, k$ and from $fj = j$ we have $j \mathcal{J} f$. Since $j$ and $k$ were arbitrary elements of $J$, it follows that $J$ is indeed a $\mathcal{J}$ class. Since $J$ contains $f$, an idempotent, it is regular. This completes the proof of items 2 and 3. Since the set of minimal functions is a regular $\mathcal{J}$ class, $(Y; J)$ is an R-representation. Thus by the Minimal Neighborhoods Lemma, minimal sets are minimal neighborhoods.

If $s \in S$ and $s$ is not constant, then there exists a minimal neighborhood $U$ contained in $s(Y)$. Let $e \in J$ be an idempotent such that $e(Y) = U$. Note that $es$ is not constant and in fact a minimal function. That is, $es \in J$. From $s \geq_J es$, it follows that $J$ is the unique cover of $C$, where $C$ is the set of the constant functions of $S$.

We prove item 5. Observe that from item 4, we have that for all $s \in S$ and all $j \in J$, $sj$ is either constant or $sj \in J$ from which it follows that $\rho$ is an S-compatible equivalence relation. Since there exist minimal sets, $\rho$ is not $\Delta$. The primitivity of $(Y; S)$ implies that $\rho = \nabla$.

We prove item 6. Since $J$ is a regular $\mathcal{J}$ class, $(Y; J)$ is by definition an R-representation; moreover, since $J$ consists of minimal functions, $(Y; J)$ is a c-representation. Finally $\rho = \nabla$ implies that every element $b \in Y$ is contained in some minimal set. That is, $(Y; J)$ is range-covered. □

For the remainder of the section, $J$ will refer to the (regular $\mathcal{J}$ class) of minimal functions of $(Y; S)$. The next group of lemmas culminates in a Proposition 4.6 which states that if $(Y; S)$ is primitive and $|Y| > 2$, then $(Y; J)$ is primitive.

**Lemma 4.4.** *Let $(Y; S)$ be a finite primitive representation and $|Y| > 2$. Then $(Y; J)$ is reduced. In fact, for $a \neq b \in Y$, there exists an idempotent $e \in J$ such that $e(a) \neq e(b)$.*

*Proof.* Since $|Y| > 2$, the set $J$ of minimal functions is non-empty. By Lemma 4.3, $J$ is a regular $\mathcal{J}$ class. Thus exists an idempotent $f$ in $J$. Let $U = f(Y)$. By definition of minimal set, $|U| \geq 2$. Let $\epsilon_U = \{(a, b) \in Y^2 : \forall s \in S\, (s(Y) \subset U \to s(a) = s(b))\}$. Observe that $\epsilon_U$ is an S-compatible equivalence relation.

For $u \neq v \in U$, we have $(u, v) \neq \epsilon_U$ so $\epsilon_U$ is not the universal relation. Primitivity of $(Y; S)$ implies $\epsilon_U = \Delta$. Thus for every pair of distinct elements $a, b \in Y$, there exists $k \in S$ satisfying $k(Y) \subseteq U$ and $k(a) \neq k(b)$. Since $J$ is a regular $\mathcal{J}$ class, there exists an idempotent $e$ such that $e \mathcal{L} k$. By Observation 2.1, $ker(e) = ker(k)$; hence, $e(a) \neq e(b)$. □

Our proof follows [6], a more general result that has played an important part in the development of the theory of finite algebras.

**Lemma 4.5.** *Suppose that $(Y; S)$ is primitive, $e$ is an idempotent in $J$ and $G = eJe \cap J$, a subgroup of $S$. Then the group representation $(e(Y); G)$ is primitive. Thus the c-representation $(Y; J)$ m-expands a primitive representation.*

*In particular, if $a \neq b \in Y$ are contained in the minimal neighborhood $e(Y)$ of $Y$, then in $(Y; S)$, the S-compatible relation $Cg(a, b)$ contains all pairs in $e(Y)^2$*

*Proof.* We show that $(e(Y); G)$ is primitive by contradiction. Assume that $\theta$ is a non-trivial G-compatible relation on $e(Y)$. In that case, we can choose three distinct elements $u, v, w \in e(Y)$ such that $(u, v) \in \theta$ and $(u, w) \notin \theta$. By Lemma 4.1, primitivity of $(Y; S)$ guarantees the existence of a sequence $u = u_1, \ldots, u_n = w$ of $Y$ such that each pair of adjacent elements is the image of $\{u, v\}$ under some $h \in S^1$. Apply $e$ to the left and right of each such $h$. The resulting sequence of elements is contained in $e(Y)$ and each pair of adjacent elements is the image of $\{u, v\}$ under an element of $eSe$. So we can intrepret the resulting sequence as one which witnesses in $(e(Y); G)$ the inclusion of $(u, w) \in \theta$, contradicting the choice of $(u, w)$. We have shown that $(e(Y); G)$ is primitive. In particular, for $a \neq b \in e(Y)$ any pair $(c, d)$ contained in $e(Y)^2$ lies in the equivalence relation generatd by $\{(g(a), g(b)) : g \in G\}$. Since $G \subseteq S$, the last sentence of the lemma follows. $\square$

**Proposition 4.6.** *If $(Y; S)$ is a faithful primitive representation, $|Y| > 2$, then $(Y; J)$ is a faithful primitive c-representation.*

*Proof.* Suppose that $(Y; S)$ is a faithful primitive representation and $|Y| > 2$. Let $a \neq b \in Y$. By Lemma 4.4, there exists an idempotent $e \in J$ such that $e(a) \neq e(b)$. Let $U = e(Y)$. By Lemma 4.5, we have that $Cg(e(a), e(b))$ contains $U^2$. In particular, $Cg(a, b)$ contains $U^2$. Let $V \in M(Y; S)$ and let $f \in J$ be such that $f(Y) = V$. Since $J$ is a $\mathcal{J}$ class, there exists $s, t \in S$ such that $f = set$. It follows that $s(U) = V$ and that $V^2$ is contained in $Cg(a, b)$. In particular, the $S$-compatible relation $\rho$ is contained in $Cg(a, b)$. By Lemma 4.3, item 5, $\rho = \nabla$. Hence, $Cg(a, b) = \nabla$, and the lemma follows. $\square$

### 4.2.1 Primitivity and translational hulls

For a primitive representation $(Y; S)$, where $|Y| > 2$ and $J$ is the set of minimal functions, we show in Propositions 4.7, 4.8 that $(Y; J)$ influences $(Y; S)$ in an interesting way.

**Lemma 4.7.** *Suppose $(Y; S)$ is a faithful primitive representation and $|Y| > 2$. For each pair $s, t \in S$ of distinct non-constant elements, there exists $j, k \in J$ such that $js \neq jt$ and $sk \neq tk$ and .*

*Proof.* Suppose $s, t \in S$ are distinct non-constants. Since they are distinct, there exists $y \in Y$ such that $s(y) \neq t(y)$. By Lemma 4.4, $(Y; J)$ is reduced; hence, there exists $j \in J$ such that $js(y) \neq jt(y)$. Thus $js \neq jt$. By Lemma 4.3, item 5, we have $(Y; J)$ is range-covered and as remarked, by regularity of $J$, the set $j(Y)$ is an image of an idempotent of $J$. Hence, for $y$ above, there exists an idempotent $k \in J$ such that $k(y) = y$. We have $sk(y) \neq tk(j)$. Thus $sk \neq tk$. $\square$

**Proposition 4.8.** *Suppose that $(Y; S)$ is a faithful primitive representation, $|Y| > 2$ and $J$ is the set of minimal functions of $(Y; S)$. Then*

1. *$S^Y$ is isomorphic to a subsemigroup of $\Omega(J^Y)$.*

2. *The set of minimal functions of the primitive faithful representation $(Y; \Omega(J^Y))$ is $J$.*

*Proof.* By Lemma 4.3, item 6, we have $J$ is the unique cover of the $\mathcal{J}$ class $C$ of constant functions of $(Y; S)$ in the partially ordered set of $\mathcal{J}$ classes. Thus if $s \in S$, then $sJ \cup Js \subseteq J \cup C$. It follows that for all $s \in S^Y$, $sJ^Y \cup J^Y s$ is contained in $J^Y$. In particular, $S^Y \leq \Omega(J^Y)$.

For the second item, note that since $\Omega(J^Y$ contains $J^Y$, and $(Y;J)$ is primitive by Lemma 4.6, we have $(Y;\Omega(J^Y))$ is primitive. Suppose that $u \in \Omega(J^Y)$ is a minimal function of $(Y;\Omega(J^Y))$. We show that $u \in J$. By Lemma 4.3, item 5, $\rho = \nabla$ from which it follows that if $u$ is non-constant, then there exists $k \in J$ such that $uk$ is not constant. But $u \in \Omega(J^Y)$ implies that $uk \in J$. By minimality of $u$, we have that $uk$ is also minimal; hence, the image of $u$ and the image of $uk$ are equal. By the regularity of $J$, there exists an idempotent $e \in J$ such that the image of $u$ is equal to the image of $e$. We have $eu = u$. Once again using the fact that $u \in \Omega(J^Y)$, we have that $u \in J$. □

**Lemma 4.9.** *Let $(Y;S)$ be a faithful primitive representation, $|Y| > 2$. Then there exists a faithful primitive group representation $(X;G)$, a reductive c-ramified Rees matrix $P$ over $(X;G)$ such that $(Y;J^Y)$ is equivalent to $(V_{\theta_P}; L((X;G);P))$.*

*Proof.* By Lemma 4.3 and Lemma 4.4, we have that $(Y;J)$ is a faithful reduced c-representation. By the Fundamental Lemma of R-representations and Lemma 3.15, there exists a faithful group representation $(X;G)$ and a reductive c-ramified matrix $P$ over $(X;G)$ such that $(Y;J)$ is equivalent to $(V_{\theta_P}; J((X;G);P))$ where $(V_{\theta_P}; L((X;G);P))$ m-expands $(X;G)$. Replacing $J$ by $J^Y$ results in a change from $J((X;G);P)$ to $L((X;G);P)$ above. By the Minimal Neighborhoods Lemma, an R-representation m-expands a unique up to equivalence group representation. By Lemma 4.5, $(Y;J^Y)$ m-expands a primitive group representation; hence, $(X;G)$ is primitive. □

Consider $(V_{\theta_P}; L((X;G);P))$, where $P$ is reductive and $(X;G)$ is primitive as in the previous lemma. We have $J((X;G);P) = J$ is the set of minimal functions of $(V_{\theta_P}; L((X;G);P))$. Recall that for $(V_{\theta_P}; L((X;G);P))$ the set of minimal neighborhoods is denoted by $\{X_i : i \in I\}$. Let $\rho'$ be defined as $Eq(\{X_i^2 : i \in I\})$. By Lemma 4.3, item 3, every $s \in S$ maps a minimal neighborhood to a singleton set or to another minimal neighborhood. Hence, $\rho'$ is a $J$-compatible equivalence relation on $Y$. Recall the definition of $Gr(P)$, the graph associated with c-ramified Rees matrix $P$.

**Lemma 4.10.** *Let $P$ be a reduced $m \times n$ c-ramified matrix over the group representation $(X;G)$ and consider the reduced c-representation $(V_{\theta_P}; L((X;G);P))$.*

*For $r \neq s \in I$, we have that $X_r \cap X_s \neq \emptyset$ if and only if $P^{(r)}, P^{(s)}$ is an edge of $Gr(P)$. In particular, $\rho' = \nabla$ if and only if $Gr(P)$ is connected. Thus if $(V_{\theta_P}; L((X;G);P))$ is primitive, then $Gr(P)$ is connected.*

*Proof.* Suppose $X_r$ and $X_s$ are minimal sets. Then $X_r \cap X_s \neq \emptyset$ if and only if there exist $x, w \in X$ such that $r_x/\theta_P = s_w/\theta_P$ if and only if $P^{(r)}, P^{(s)}$ is an edge of $Gr(P)$. Item 2 follows by inspection. □

**Lemma 4.11.** *Suppose that $(Y;S)$ is a faithful primitive representation and $|Y| > 2$. Then $(Y;J^Y)$ is equivalent to $(V_{\theta_P}; L((X;G);P))$ where*

1. $(X; G)$ is a primitive group representation.

2. $P$ is a reductive c-ramified Rees matrix.

3. $Gr(P)$ is connected.

*Conversely, for a group representation $(X; G)$ and a c-ramified matrix $P$ over $(X; G)$, we have $(V_{\theta_P}; L((X; G); P))$ is a faithful primitive representation if and only if items 1-3 above are satisfied.*

*Proof.* The forward direction follows from Lemma 4.9 and Lemma 4.10.

Conversely, suppose that $(V_{\theta_P}; L((X; G); P))$ satisfies items 1-3 above. Let $a \neq b \in Y$. Since $P$ is a reductive, by Lemma 3.15, $(V_{\theta_P}; L((X; G); P))$ is reduced. Hence by Lemma 4.4, there exists an idempotent $e \in J((X; G); P)$ such that $e(a) \neq e(b)$. Since $(X; G)$ is primitive, by the Minimal Neighborhoods Lemma we have that $(e(Y); eJe \cap J)$ is primitive. Since $Gr(P)$ is connected, by Lemma 4.10 $\rho' = \nabla$. Arguing as in the proof of Lemma 4.6, we have that $(V_{\theta_P}; L((X; G); P))$ is primitive. □

We prove Theorem 1.8.

*Proof.* (Proof of Theorem 1.8) Suppose that $(Y; S)$ is a faithful primitive representation, $|Y| > 2$ and $J$ is the set of minimal functions of $(Y; S)$. By Lemma 4.8, we have $S^Y$ is isomorphic to a $J$-containing subsemigroup of $\Omega(J^Y)$. By Proposition 4.11, $(Y; J^Y)$ is equivalent to a reduced c-representation $(V_{\theta_P}; L((X; G); P))$ where $(X; G)$ and $P$ satisfy items 1-3 of Lemma 4.11. By Lemma 4.8, $J$ is the set of minimal functions of $\Omega(J^Y)$. Finally by Lemma 4.7, for any $s, t \in T$, there exist $j, k \in J$ such that $sj \neq tj$ and $ks \neq kt$.

Conversely if $T$ is any $J$ containing subsemigroup of $\Omega(L((X; G); P)$ where $(X; G)$ and $P$ satisfy items 1-3, then by the converse part of Lemma 4.11, we have $(V_{\theta_P}; L((X; G); P))$ is faithful and primitive. Thus if $T$ is any subsemigroup of $\Omega(L((X; G); P)$ which contains $J((X; G); P)$, then $(V_{\theta_P}; T)$ is faithful (since $\Omega(J^Y)$ is defined as a semigroup of transformations) and primitive. □

## 5 Equivalences of primitive representations

Isomorphism classes of Rees matrix semigroups admit a sharp description. See Theorem 5.2 below. In this section we generalize Theorem 5.2 to certain transformation semigroups, the purpose being to further understand the interplay between c-ramified matrices and primitive representations.

Recall that if $G$ is group and $n$ is a positive integer, then an $n \times n$ matrix $W$ is *permutational* if $W$ is the result of taking an $n \times n$ permutation matrix $\Pi$ and replacing each 1 entry in $\Pi$ by an element of $G$. Every permutational matrix $W$ admits two useful factorizations. We have $W = D\Pi$ where $D$ is a diagonal matrix with diagonal entries in $G$ and $\Pi$ is a permutation matrix and $W = \Pi'D'$, where $\Pi'$ is a permutation matrix and $D'$ is a diagonal matrix with diagonal entries in $G$. Conversely, any matrix which admits a factorization of either type above is permutational.

**Definition 5.1.** *Let $G$ and $H$ be groups, $Q$ a regular matrix $m \times n$ over $G$, $Q'$ a regular matrix $m \times n$ over $H$. We say that $Q$ and $Q'$ are **equivalent** ($Q \equiv Q'$) if there exist $U$, an $m \times m$ permutational matrix over $G$, $V$, an $n \times n$ permutational matrix over $G$ and a group isomorphism $\phi : G \to H$ such that $Q' = \phi(UQV)$.*

Here is the theorem from the theory of (finite) Rees matrix theory which characterizes isomorphism classes of regular Rees matrix semigroups. Ther reader may want to compare Proposition 1.9 and its terminology with the well-known theorem that follows.

**Theorem 5.2.** *Let $S = \mathcal{M}(Q, G)$ and $S' = \mathcal{M}(Q', G')$ be two finite regular Rees matrix semigroups. Then $S$ is isormorphic to $S'$ if and only if $Q \equiv Q'$.*

The next proposition is stronger than Proposition 1.9. We use the remainder of the subsection to prove it.

**Proposition 5.3.** *Suppose $P$ is a c-ramified reductive matrix over $(X; G)$ and $P'$ is a c-ramified matrix over $(X'; G')$. Then the following are equivalent.*

1. *$P$ is equivalent to $P'$*

2. *$U((X;G); P)$ is isomorphic to $U((X'; G'); P')$*

3. *$L((X;G); P)$ is isomorphic to $L((X'; G'); P')$*

4. *$(V_{\theta_P}; L((X; G); P))$ is equivalent to $(V_{\theta_{P'}}; L((X'; G'); P'))$*

**Definition 5.4.** *Let $P$ be matrix over $(X; G)$. If $(Y; J)$ is an R-representation, $J \mathcal{P}$ is a Rees generating set for $J$, and $P$ is the ramified Rees matrix associated with $\mathcal{P}$, then we then say that $P$ is **definable** by $(Y; J)$.*

**Lemma 5.5.** *If $(Y; J)$ and $(Y'; J')$ are equivalent R-representations, then for any matrix $P$ which is definable by $(Y; J)$, there exists a $(Y'; J')$-definable matrix $P'$ such that $P \equiv P'$.*

*Proof.* Let $P$ be definable by $(Y; J)$ and let $\mathcal{P}$ be the ordered Rees generating set over the idempotent $e \in J$ which gives rise to the ramified matrix $P$. Let $\phi : (Y, J) \to (Y'; J')$ be an equivalence. The restriction of $\phi$ to $\mathcal{P}$ maps $\mathcal{P}$ to a Rees generating set for $(Y'; J')$ which we call $\mathcal{P}'$. Observe that the ramified matrix associated with $\mathcal{P}'$ is $\phi(P)$. By definition of equivalence for ramified matrices, we have that $\phi(P)$ is equivalent to $P$. □

For an R-representation $(Y; J)$, we provide a description of the collection of Rees generating sets for $J$ and use it to to describe the $(Y; J)$-definable matrices. Suppose $\mathcal{P}$ is based at the idempotent $e \in J$ and $\mathcal{P} = \mathcal{I} \cup \mathcal{K} \cup H$ where $\mathcal{I} = \{f_r : r \in I\}$, $\mathcal{K} = \{g_\lambda : \lambda \in \Lambda\}$ and $H = eJe \cap J$. Let $e'$ be an idempotent in $J$. Since $e \mathcal{J} e'$, there exist $s, t \in J$ such that $set = e'$. Let $\{b_r \in J : r \in I\}$ and $\{c_\lambda \in J : \lambda \in \Lambda\}$ be indexed subsets of $J$ satisfying the following:

- For each $r \in I$, $b_r \,\mathcal{H}\, e$
- For each $\lambda \in \Lambda$, $c_\lambda \,\mathcal{H}\, e$.

Using Green's Lemma, it follows that

- $\mathcal{I}' = \{f_r b_r t : r \in I\}$ contains a representative from each $\mathcal{H}$ class of $L_{e'}$
- $\mathcal{K}' = \{sc_\lambda g_\lambda : \lambda \in \Lambda\}$ contains a representative from each $\mathcal{H}$ class of $R_{e'}$.

By Lemma 2.4, $sHt$ is a maximal subgroup of $J$; hence, $\mathcal{I}' \cup \mathcal{K}' \cup sHt$ is a Rees generating set. Moreover since for all $r \in I$ and all $\lambda \in \Lambda$ we have let $b_r$ and $c_\lambda$ range over $eJe \cap J$, by Green's Lemma again, **all** ordered Rees generating sets are described by this method. That is, if $\mathcal{P}$, $\mathcal{P}'$ are ordered Rees generating sets for $J$ over idempotents $e$ and $e' = set$ respectively, then "transforming" $\mathcal{P}$ into $\mathcal{P}'$ involves (independently) choosing sequences in $H_e$, $\{b_r : r \in I\}$ and $\{c_\lambda : \lambda \in \Lambda\}$, and a permutation $\Pi$ of $I$, a permutation $\Phi$ of $\Lambda$ ( the permutations come into the description because we are dealing with ordered Rees generating sets). The next lemma is a consequence of the the claims of this paragraph, along with the already noted observation that the product of a permutation matrix and a diagonal matrix is a permutational matrix.

**Proposition 5.6.** *Let $(Y; J)$ be an R-representation with ordered Rees generating set $\mathcal{P}$ over the idempotent $e$, whose associated ramified ramified matrix is $P$ (so $P$ is definable by $(Y; J)$). Then $P'$ is definable by $(Y; J)$ if and only if there exist*

- *An idempotents $e' \in J$, elements $s, t \in J$ such that $set = e'$*
- *$U$, an $m \times m$ permutational matrix over $(X; G)$, $V$, an $n \times n$ permutational matrix over $(X; G)$*

*such that $P' = \beta(UPV)$ where $\beta : (e(Y); eJe \cap J) \to (e'(Y); e'Je' \cap J)$ is the equivalence determined by $s, t$ (described in the Minimal Neighborhoods Lemma) given by $e(y) \to se(y)$ and $eje \to sejet$. In particular, if $P'$ is definable by $(Y; J)$, then $P'$ is equivalent to $P$.*

*Proof.* (Proof of Proposition 5.3) Items 3 and 4 are equivalent in view of the fact that isomorphism between the semigroups in question induces an obvious equivalence among the representations and conversely.

Suppose that item 3 holds. Since the isomorphism of $L((X; G); P)$ and $L((X'; G'); P')$ induces an equivalence between $(V_{\theta_P}; L((X; G); P))$ and $(V_{\theta_{P'}}; L((X'; G'); P'))$, by Lemma 5.5, we have $P \equiv P'$. This proves item 3 implies item 1.

We show that item 2 implies item 3. Suppose $U((X; G); P)$ is isomorphic to $U((X'; G'); P')$ under an isomorphism $\beta$. It is easy to verify that $\beta$ maps the pairs of the maximal deflation $\theta_P$ to the pairs of the maximal deflation $\theta_{P'}$. Thus $\beta$ induces an isomorphism from $L((X; G); P)$ to $L(X'; G'); P')$.

We prove item 1 implies item 2. Let $P$ be equivalent to $P'$. Assume first that $(X; G) = (X'; G')$ and that the equivalence is of the form $P' = UPV$ (i.e. there is no group representation equivalence $\beta$ involved).

Denote the non-zero entry of the $i$th column of $U$ by $u_i$ and denote the non-zero entry $\lambda$th row by $v_\lambda$. We define a permutation $\Pi$ by on $I$ by the rule $i \to \Pi(i)$ if $u_{\Pi(i),i}$ is non-zero entry of the i column of $U$ (thus $u_{\Pi(i),i} = u_i$); similarly, we define $\Phi$, a permutation on $\Lambda$, $\lambda \to \Phi(\lambda)$ if $v_{\lambda,\Phi(\lambda)}$ is the non-zero entry of the $\lambda$ row of $V$. We wil make use of the following observation later. For all $i \in I$ and $\lambda \in \Lambda$, observe that

$$p_{\lambda,i} = (v_{\Phi^{-1}(\lambda)})^{-1} q_{\Phi^{-1}(\lambda),\Pi^{-1}(i)} (u_{\Pi^{-1}(i)})^{-1}$$

We define a map $\phi$ from $V = V(m, X)$ to $V' = V(m, X')$ by the rule $i_x \to \Pi^{-1}(i)_{(u_{\Pi^{-1}(i)})^{-1}(x)}$. Clearly $\phi$ is a bijection from $V$ to $V'$. We extend the domain of $\phi$ to $J((X;G);P)$ by sending

$$[r, g, \lambda] \to [\Pi^{-1}(r), (u_{\Pi^{-1}(r)})^{-1}(g)(v_{\Phi^{-1}(\lambda)})^{-1}, \Phi^{-1}(\lambda)]$$

Using the displayed equation above, it is not difficult to check that the mapping $\phi$ (part of which is displayed directly above) is indeed an isomorphism from $U((X;G);P)$ to $U((X';G');P')$.

Thus for the case $P' = UPV$ we have that $U((X;G);P)$ is isomorphic to $U((X;G);P')$. Observe that for any equivalence $\beta : (X;G) \to (X';G')$, if $P' = \beta(P)$, then $U((X;G);P)$ is isomorphic to $U((X';G');P')$ via the map $(r, g, \lambda) \to (r, \beta(g), \lambda)$ and $i_x \to i_{\beta(x)}$. Using the transitivity of the isomorphism relation on semigroups, we have shown that item 2 implies item 1. □

## 6 R-representations, transitive and cyclic c-representations

In this section by defining a variation of the ramified Rees representation, we provide complete the description of the faithful R-representations. We then describe the transitive and cyclic c-representations.

Let $J^0$ be completely 0-simple and let $(Y; J)$ be a faithful representation of $J$. As usual we index the $\mathcal{L}$ classes of $J$ by the set $\Lambda = \{1, \ldots, \lambda, \ldots, m\}$. Since $(Y; J; \gamma)$ is faithful, by Lemma 3.7, for any pairs of elements of $j, k \in J$ $ker(j) = ker(k)$ if and only if $j \mathcal{L} k$. We let $\pi_1, \ldots, \pi_m$ be an enumeration of the set of kernels of $\gamma(J)$.

The next observation reduces the complicatedness of establishing that two representations are equivalent in the case that one of the representations is faithful.

**Lemma 6.1.** *Suppose $(Y; S)$ and $(Z; T)$ representations. If either $(Y; S)$ or $(Z; T)$ is faithful,then a bijection $\alpha : Y \cup S \mapsto Z \cup T$ is an equivalence if for all $y \in Y$ and all $s \in S$, $\alpha(s(y)) = \alpha(s)\ (\alpha(y))$.*

*Proof.* Assume that $(Y; S)$ and $(Z; T)$ are representations and $\alpha$ is a bijection satisfying for all $y \in Y$ and all $s \in S$, $\alpha(s(y)) = \alpha(s)\ (\alpha(y))$. We leave it to the reader to check that $(Y; S)$ is faithful if and only if $(Z; T)$ is faithful.

So assume $(Z; T)$ is faithful and that for all $y \in Y$ and all $s \in S$, $\alpha(s(y)) = \alpha(s)\ (\alpha(y))$. To see that $u, v \in S$ and $uv \in S$ imply $\alpha(uv) = \alpha(u)\ \alpha(v)$, check that the hypotheses guarantee that for all $z \in Z$, $\alpha(uv)(z) = \alpha(u)\ \alpha(v)(z)$. Since $(Z;T)$ is faithful, it follows that $\alpha(uv) = \alpha(u)\ \alpha(v)$. $\square$

**Definition 6.2.** *Suppose $(Y; J; \gamma)$ is a faithful R-representation such that $J$ has $m$ $\mathcal{L}$ classes indexed by $\{1, \ldots, \lambda, \ldots, m\}$. We define a new R-representation of $J$, $(Y, J; \gamma)^{ker} = (Y^m; J; \gamma^{ker})$ as follows. Suppose $j \in J$ and that $j$ is in the $\lambda$th $\mathcal{L}$ class of $J$. For $(a_1, \ldots, a_m) \in Y^m$, let $\gamma^{ker}(j)(a_1, \ldots, a_m) = (\gamma(j)(a_\lambda), \ldots, \gamma(j)(a_\lambda))$.*

**Lemma 6.3.** *For any R-representation $(Y; J)$, $(Y; J)^{ker}$ is an R-representation.*

*Proof.* Suppose that $j, k \in J$, $jk \neq 0$, and $k$ is in the $\lambda$ $\mathcal{R}$ class of $J$. Then by item 1 of Green's Lemma, $jk$ is in the same $\mathcal{L}$ class as $k$. Thus for $(a_1, \ldots, a_m) \in Y^m$, we have $\gamma^{ker}(jk)((a_1, \ldots, a_m) = (\gamma(jk)(a_\lambda), \ldots, \gamma(jk)(a_\lambda))$. On the other hand, $\gamma^{ker}(j)\gamma^{ker}(k)(a_1, \ldots, a_n) = \gamma^{ker}(j)((\gamma(k)(a_\lambda), \ldots, \gamma(k)(a_\lambda)) = (\gamma(j)\gamma(k)(a_\lambda), \ldots, \gamma(j)\gamma(k)(a_\lambda))$ $= (\gamma(jk)(a_\lambda), \ldots, \gamma(jk)(a_\lambda))$, as desired. $\square$

**Definition 6.4.** *The set $D = \{(a, \ldots, a) : a \in Y\}$ is referred as the set of **diagonal elements** of $Y^m$. Note that for any $j \in J$, $\gamma^{ker}(j)$ maps $Y^m$ into $D$. A subrepresentation of $(Y; J)^{ker}$ which contains $D$, the diagonal elements of $Y^m$ will be called a **diagonal subrepresentation** of $(Y; J)^{ker}$.*

We will show that every R-representation can be approximated by a $(Y; J)^{ker}$ representation, where $(Y; J)$ is a ramified Rees representation.

**Definition 6.5.** *Let $(Y; S)$ be any representation. A representation of $(Z; T)$ is deflation-equivalent to (Y;S) if there is a finite sequence of representations starting with $(Y; S)$ terminating in $(Z; T)$ such that for each adjacent pair in the sequence, either the pair is equivalent, or one of the pair is a deflation of the other.*

**Proposition 6.6.** *Let $J^0$ be a completely 0-simple semigroup and let $(Y; J; \gamma)$ be a faithful representation of $J$. Then there exists a faithful group representation $(X; G)$, a ramified $m \times n$ matrix over $(X; G)$ such that $(Y; J; \gamma)$ is deflation-equivalent to a diagonal subrepresentation of $(V_{\theta_P}; J((X; G); P))^{ker}$.*

*Proof.* Let $(Y; J; \gamma)$ be an R-representation. Since the proposition we are proving is a statement about R-representations up to deflation-equivalence we can assume without loss of generality that $(Y; J; \gamma)$ is reduced.

Condider the subrepresentation $(J(Y); J; \gamma|_{J(Y)})$. Note that since $J$ is regular, we have $J^2 = J$, from which it follows that $(J(Y); J; \gamma|_{J(Y)})$ is range-covered. Both $(J(Y); J; \gamma|_{J(Y)})$ and $(Y; J; \gamma)$ are m-expansions (up to equivalence) of $(X; G)$. Since $(Y; J; \gamma)$ is faithful, so is $(X; G)$. Thus $(J(Y); J; \gamma|_{J(Y)})$ is an m-faithful, range-covered, reduced R-representation. By the Fundamental Lemma of R-representations, $(J(Y); J; \gamma|_{J(Y)})$ is equivalent, by an equivalence $\beta$, to $(V_{\theta_P}; L((X; G); P))$, where $P$ is a ramified matrix over the group representation $(X; G)$.

We map $Y$ into $V_{\theta_P}^m$ by a map which we call $\kappa$. Suppose first that $y \in J(Y)$. We let $\kappa(y) = (\beta(y), \ldots, \beta(y))$, an element in $D$, the diagonal. For $y \in Y \setminus$ $\{J(Y)\}$, let $\kappa(y) = (\beta(y_1), \ldots, \beta(y_m))$, where for each $\lambda \in \Lambda$, $(y_\lambda, y) \in \Pi_\lambda$. We claim elements of the form $y_\lambda$ satifying $(y_\lambda, y) \in \Pi_\lambda$ do in fact exist. To see this, observe that $|Y/\Pi_\lambda|$ is the number of elements in a minimal neighborhood of $(Y; J; \gamma)$, $|J(Y)/\Pi_\lambda|_{J(Y)}|$ is the number of elements of a minimal neighborhood of $(J(Y); J; \gamma|_{J(Y)})$, and both of these representations m-expand $(X; G)$. Thus, by the Minimal Neighborhoods Lemma, $|Y/\Pi_\lambda| = |X| = |J(Y)/\Pi_\lambda|_{J(Y)}|$. The claim follows.

Since $(Y; J; \gamma)$ is by assumption reduced, it follows that $\kappa$ is one-to-one on $Y$. Note that $\kappa(J(Y)) = D$. We extend the domain of $\kappa$ to $J$. For $j \in J$, if $j$ is in the $\alpha$-th $\mathcal{L}$ class of $J$, let $\kappa(j)(a_1, \ldots, a_m) = (\beta(j)(a_\alpha), \ldots, \beta(j)(a_\alpha))$.

We show that $\kappa : (Y; J) \to (V_{\theta_P}; J((X; G); P))^{ker}$ is an equivalence into a diagonal subrepresentation of $(V_{\theta_P}; J((X; G); P))^{ker}$. Since $(Y; J; \lambda)$ is faithful, it suffices to show that for all $j \in J$ and $y \in Y$, we have $\kappa(j(y)) = \kappa(j)\kappa(y)$. Since $j(y) \in J(Y)$, we have $\kappa(j(y)) = (\beta(j(y)), \ldots, \beta(j(y)))$. On the other hand, $\kappa(j)\kappa(y) = \kappa(j)(\beta(y_1), \ldots, \beta(y_m))$ where $(y_\alpha, y) \in \Pi_\alpha$ $\alpha \in \Lambda)$. We have $\kappa(j)(\beta(y_1), \ldots, \beta(y_m)) = (\beta(j)(\beta(y_\alpha)), \ldots, \beta(j)(\beta(y_\alpha))) = (\beta(j(y_\alpha)), \ldots, \beta(j(y_\alpha)))$ (since $\beta$ is an equivalence) $= (\beta(j(y)), \ldots, \beta(j(y)))$, the last equality holding since $j(y) = j(y_\alpha)$. We have shown that $\kappa$ is an equivalence thereby completing the proof of the lemma. □

We make our characterization of faithful R-representations exact. For an R-representation $(Y; J)$ and a set $Z$ containing a distinguished element $\epsilon \in Z$, we let $(Y; J; \gamma^{ker,Z})^{ker,Z}$ be the action of $J$ on the set $J(Y)^m \times Z$ where for all $j \in J$, $(u, z) \in J(Y)^m \times Z$, $\gamma^{ker,Z}(u, z) = (\gamma^{ker}(u), \epsilon)$. It is straightforward to verify that $(Y; J; \gamma^{ker,Z})^{ker,Z}$ is in fact an R-representation. A diagonal subrepresentation of $(Y; J; \gamma^{ker,Z})^{ker,Z}$ is a subrepresentation containing $\{(a, \ldots, a, \epsilon) : a \in Y\}$ (i.e. $D \times \epsilon$).

**Proposition 6.7.** *Let $(Y; J; \gamma)$ be a faithful R-representation. Then there exists $(X; G)$ a faithful group representation, $P$ a ramified matrix over $(X; G)$, $\alpha \leq \theta_P$, a set $Z$ containing a distinguished element $\epsilon$ such that $(Y; J; \gamma)$ is equivalent to a diagonal subrepresentation of $(V_\alpha; J((X; G); P))^{ker,Z}$.*

*Proof.* The construction is similar to the one in Proposition 6.6. We map $\gamma(J)(Y)$ to the "diagonal" elements of $\{(a, \ldots, a, \epsilon) : a \in V_\alpha$, having made use of the fact that $(J(Y); J; \gamma|_{J(Y)})$ is equivalent to a ramified Rees representation $(V_\alpha; J((X; G); P))$. We mimic the proof of Proposition 6.6 further by initially sending $y \to \kappa(y) = (a_1, \ldots, a_m)$, where $(y, a_\lambda) \in Pi_\lambda$ $(\lambda \in \Lambda)$. Since we are no longer assuming that $(Y; J; \lambda)$ is reductive, so $\kappa$ is not necessarily one-to-one. To remedy this, we define a large enough set $Z$ with distinguished element $\epsilon$, and an extension of $\kappa$, which we call $\kappa^Z$ so that $\kappa^Z(Y) \to V_\alpha^m \times Z$ is one-to-one, where we define $\gamma^{ker,Z}(a_1, \ldots, a_m, z) = (\gamma(a_1), \ldots, \gamma(a_m), \epsilon)$. The remainder of the proof, arguing that we have constructed is indeed an equivalence, closely follows the proof of Proposition 6.6 and is left to the reader. □

**Remarks 6.8.** *In light of Lemma 6.7, each R-representation $(Y; J)$ is equivalent to we will call a **two-sorted Rees matrix object**: The transformations* are the of the form $M(j) \in J((X;G);P)$ where $J((X;G);P)$ are the $m \times n$ monomial matrices over a group representation $(X;G)$; the carrier set is approximated by a set of matrices, the $n \times n$ column monomial matrices over $V_\alpha$. We add a place holder set $Z$ with distinguished element $\epsilon$ so that $y \in Y$ is mapped to $(A(y), z(y))$ where $A(y)$ is an $n \times n$ column monomial matrix, $z(y) \in Z$ and for $j$ in $J$, $M(j)(A(y), z(y)) = (MPA(y), \epsilon)$.

## 6.1 Some applications involving transitive and cyclic c-representations

Cyclic and transitive representations of finite simple semigroups were treated in the some early papers of the algebraic theory of semigroups ( see [10] and [11]). Of course every representation of a simple semigroup is vacuously a c-representation. Every transitive representation is a range-covered representation so the Fundamental Lemma of R-representations give us some headway with the transitive R-representations. We give a description of the transitive and cyclic c-representations.

### 6.1.1 Transitive and cyclic c-representations

Let $(X;G)$ be a group representation. Recall that for $x \in X$ the $G$-transitivity class $G_x$ of $x$ under the representation $(X;G)$ is $\{g(x) : g \in G\}$. Suppose $P$ is a $m \times n$ c-ramified matrix over $(X;G)$. For $r \in I$, let $K_r$ denote the constants in $P^{(r)}$, the r column of $P$ and let $G_{x,r}$ be the union of the $G$-transitivity classes of $K_r \cup \{x\}$.

**Proposition 6.9.** *Let $\mathcal{M}\ (G,Q)$ be a regular Rees matrix semigroup and let $J$ denote its maximal $\mathcal{J}$ class. Suppose $(Y;J)$ is a c-representation of $J$. Then $(Y;J)$ is transitive if and only if $(Y;J)$ is equivalent to a c-ramified Rees representation $(V_\alpha; J((X;G);P))$ where $P$ is a c-ramification of $Q$ over $(X;G)$ and for each $x \in X$ and each $s \in I$, the set $\{r_y/\alpha : r \in I, y \in G_{x,s}\} = V_\alpha$. In particular, if $\mathcal{M}\ (G;Q)$ is simple (so $J = \mathcal{M}\ (G;Q)$), then $(Y;J)$ is transitive if and only if $(Y;J)$ is equivalent to a c-ramified Rees representation $(V_\alpha; M((X;G);P))$ where $P = \gamma(Q)$ and $\{r_y/\alpha : r \in I, y \in G_x\} = V_\alpha$.*

*Proof.* We have that $(V_\alpha; J((X;G);P)$ is transitive if and only if for all $a_x, s_y \in V$, there exists $[r, g, \lambda]$ such that $r_{gP(\lambda,a)(x)}/\alpha = [r, g, \lambda](a_x)/\alpha = s_y/\alpha$. The proposition now follows by inspection. □

**Remarks 6.10.** *Given as a starting point a Rees matrix semigroup $\mathcal{M}\ (G,Q)$ and a representation $(X;G)$ of its group $G$, it would be interesting to be able to find an effective way to construct (all) faithful transitive c-representations which m-expand $(X;G)$. In many cases, the proposition above does lend itself to quickly producing faithful transitive representations if these exist, but in the general case, as far as we know, a procedure that produces all the faithful transitive c-ramifications of $Q$ over $(X;G)$ that satisfy the conditions of the proposition, can not be done via any type of implementable algorithm.*

**Problem 6.11.** *Consider the following computational complexity problem. Let $(X;G)$ be a faithful group representation. An instance is an abstract Rees matrix semigroup with structure matrix $G$, $\{\mathcal{M}(G,Q)\}$. The question is whether there exists a c-ramification $P$ of $Q$ over $(X;G)$ such that $(V_{\theta_P}; J((X;G);P))$ is faithful and transitive. So $\mathcal{M}(G,Q)$ will range over the Rees matrix semigroups with structure group $G$; the size of the instance is is the number of elements of $\mathcal{M}(G,Q)$). Note that the problem is in NP since given $P$, a c-ramification of $Q$ over $(X;G)$, using Proposition 6.9, we can run a polynomial time check to determine if $(V_{\theta_P}; J((X;G);P))$ is transitive. Moreover, $(V_{\theta_P}; J((X;G);P))$ is faithful if and only if $P$ is reductive and this too can be quickly checked.*

*Does there exist a faithful group representation such that the above complexity problem is NP-complete?*

We treat cyclic representations. There are two cases since cyclic c-representations are not necessarily range-covered. To treat range-covered cyclic c-representations, we need change only the quantifiers in Proposition 6.9. For the non-range covered cyclic c-representations, we define an expansion of the usual ramified Rees c-representation $(V_\alpha; J((X;G);P))$. Let $i = ((r_1)_{x_1}, \ldots, (r_m)_{x_m})$ be a m-tuple with entries in $V_\alpha$. We extend the domain of an arbitrary $J \in M((X;G);P)$ to $V_\alpha \cup \{i\}$ as follows. For $[r,g,\lambda] \in J((X;G);P))$, let $[r,g,\lambda](i) = r_y/\alpha$ where $y = gP(\lambda, r_\lambda)(x_\lambda)$. We denote $V_\alpha \cup i$ by $V_\alpha^i$. Note that $(V_\alpha^i; J((X;G);P))$ is a c-representation (in fact, it is equivalent to a diagonal subrepresentation of $(V_\alpha; J((X;G);P))^{ker}$). A representation of the form $(V_\alpha^i; M((X;G);P))$ will be called a **initial state ramified Rees representation** or just an **i-ramified Rees representation**. Let $G_i$ be the union of the $G_{x_\lambda, r_\lambda}$, as $\lambda$ varies from 1 to $m$. The proof of the following proposition follows closely along the lines of the proof of Proposition 6.9 and is left to the reader.

**Proposition 6.12.** *Let $(Y;J)$ be a c-representation of $J$. Then $(Y;J)$ is cyclic if and only if it is equivalent to one of the following two types of c-representation.*

1. *a c-ramified Rees representation $(V_\alpha; J((X;G);P))$ where $P$ is a c-ramification of $Q$ over $(X;G)$ and there exists a column $P^{(r)}$ of $P$ ($r \in I$) and an $x \in X$, such that $\{r_y/\alpha : r \in I, y \in G_{x,r}\} = V_\alpha$ (in which case, $r_x$ is a cyclic generator).*

2. *An i-ramified Rees c-representation such that $i = ((r_1)_{x_1}, \ldots, (r_m)_{x_m})$ where the set $\{r_y/\alpha : r \in I, y \in G_i\}$ (in which case, $i$ is a cyclic generator).*

# 7 Conclusion, Problems

We have been able to uncover part of the underlying combinatorial structure of primitive representations, the c-ramified matrix invariant. Representations can quickly be be constructed from c-ramified matrices; moreover, these representations, the c-ramified Rees representations, are rather familiar algebraic objects, matrices acting on a set of vectors. We have concentrated on primitive representations of the general interest such representations hold; however, our main theorem, Theorem 1.8, minus the item that requires that the group representation $(X; G)$ be primitive, actually gives a complete description of the tame representations. We refer the reader to [3] and [2] for further information concerning tame representations and tame congruence theory. In the next paper,[7], in this series, the first author provides a description a description of primitive representations as matrix actions on sets of vectors, extending the matrix description of the action beyond the minimal functions. We say that $S$ is a **primitive semigroup** if $S$ has a faithful primitive representation. We describe the finite primitive semigroups in [9] and explore some computational complexity issues related to primitive representations ([9]). We note that our characterization of the finite primitive representations constitutes a characterization of clones of the finite simple unary algebras. In [8] we explore how the c-ramified matrix invariant of primitive representations interact with tame congruence to influence the structure of finite simple algebras.

We have described the finite primitive representations; primitive representations on infinite sets is a problem of enormously larger scope. Certain of the results here do generalize to classes of infinite primitive representations. With the following two problems, we hope to begin the search for interesting general results concerning infinite primitive representations.

The next problem concerns the possibility of locating one boundary line for the infinite primitive representations.

**Problem 7.1.** *Consider representations $(Y; S)$ of the following type. For every $s \in S$, there exists a positive integer $n$ such that $s^n$ is constant ($n$ may depend on $s$) Call such a representation a **nil-by-constant** representation. Do there exist primitive faithful nil by constant representations $(Y; S)$ where $|Y|$ is infinite.*

Note using Theorem 1.8, it is not difficult to show that if $Y$ is finite and $|Y > 2$, then there do not exist examples of faithful primitive nil by constant representations. We suspect that the there are no infinite examples either. Obviously the semigroup part of a nil-by-constant representation is group-free.

**Problem 7.2.** *What can be said about infinite faithful primitive representations $(Y; S)$ where for some positive integer $n$, $S$ is contained in the variety generated by $T_n$ and $S$ is group-free?*

Perhaps the primitive representations of this type are classifiable.